\documentclass[12pt]{article}
\usepackage[explicit]{titlesec}
\usepackage[none]{hyphenat}
\usepackage{amsmath,amssymb,amsfonts,mathrsfs}

\usepackage{soul}
\setul{.6pt}{.4pt}

\usepackage{geometry}
\usepackage{fancyhdr}
\usepackage[labelfont={footnotesize,bf} , textfont=footnotesize]{caption}
\titleformat{\section}
{\normalfont\bfseries}{\thesection.}{0.5em}{\MakeUppercase{\textbf{#1}}}
\titlespacing\section{0pt}{0pt}{-10pt}
\titleformat{\subsection}
{\normalfont\itshape}{\thesubsection.}{0.5em}{\textit{#1}}
\titlespacing\subsection{0pt}{0pt}{-8pt}

\makeatletter
\newcommand\sixteen{\@setfontsize\sixteen{17pt}{6}}
\renewcommand{\maketitle}{\bgroup\setlength{\parindent}{0pt}
	\begin{flushleft}
		\sixteen\bfseries \@title
		\medskip
	\end{flushleft}
	\textit{\@author}
	\egroup}
\makeatother

\usepackage{caption}
\usepackage{graphicx,subfigure}
\graphicspath{{figures/}}

\usepackage{booktabs}
\usepackage{multirow}
\usepackage{array}

\usepackage[american]{babel}
\usepackage{microtype} 
\usepackage{enumerate}

\usepackage{listings}
\usepackage{color}
\definecolor{dkgreen}{rgb}{0,0.6,0}
\definecolor{gray}{rgb}{0.5,0.5,0.5}
\definecolor{mauve}{rgb}{0.58,0,0.82}
\title{A multi-grid sampling multi-scale method for crack initiation and propagation}

\author{
	Zhenxing Cheng$^{a}$, Hu Wang$^{a*}$ \medskip \\
	$^{a}$State Key Laboratory of Advanced Design and Manufacturing for Vehicle Body, Hunan University, Changsha, 410082, PR China \\}

\begin{document}
	
	\vspace*{.01 in}
	\maketitle
	\vspace{.12 in}
	
	\section*{abstract}
	
	In this study, a multi-grid sampling multi-scale (MGSMS) method is proposed by coupling with finite element (FEM), extended finite element (XFEM) and molecular dynamics (MD) methods.Crack is studied comprehensively from microscopic initiations to macroscopic propagation by MGSMS method. In order to establish the coupling relationship between macroscopic and microscopic model, multi-grid FEM is used to transmit the macroscopic displacement boundary conditions to the atomic model and the multi-grid XFEM is used to feedback the microscopic crack initiations to the macroscopic model. Moreover, an image recognition based crack extracting method is proposed to extract the crack coordinate from the MD result files of efficiently and the Latin hypercube sampling method is used to reduce the computational cost of MD. Numerical results show that MGSMS method can be used to calculate micro-crack initiations and transmit it to the macro-crack model. The crack initiation and propagation simulation of plate under mode I loading is completed.

	\textit{Keywords}: Molecular dynamics, Crack propagation, Poly-crystals, Multi-scale, Extended finite element method
	
	\vspace{.12 in}
	
	
	\section{introduction}
	
	Fracture of materials usually comes from the defects of its internal micro-structure, which is caused by the coupling action of multiple failure mechanisms at macroscopic, mesoscopic and even microscopic scales. Therefore, it is necessary to understand the crack propagation under different length scales, e.g., macro-scale, micro-scale, nano-scale (or atomic scale) \cite{benz2015reconsiderations,li2015fatigue,leung2014atomistic,zhang2017mechanisms,yasbolaghi2020micro}. For the research on macro crack propagation, many numerical methods have been developed to study the crack propagation, including the finite element method (FEM) \cite{branco2015review}, boundary element method (BEM) \cite{santana2016dual}, meshless method \cite{khosravifard2017accurate}, extended finite element method (X-FEM) \cite{belytschko2009review,cheng2019exact}. For micro crack propagation, crystal plasticity finite element method (CPFEM) is a popular method to analyze the plastic deformation at the crack tip in poly-crystals \cite{zhang2018quantitative,proudhon20163d}. In addition, molecular dynamics (MD) method has also been widely used to study behaviors of micro crack during the last decades \cite{cui2014molecular,wu2015molecular,sung2015studies,chandra2016molecular,feng2018twin,chowdhury2019effects,lu2020cohesive}.
	
	In the macroscopic and mesoscopic scale, numerical methods such as the finite element method can be used to analysis the crack propagation, but it is difficult to consider the characteristics of the material micro-structure. In terms of micro scale, although micro numerical calculation methods such as molecular dynamics can be used to analysis the evolution of material micro-structure, but the computational size  is usually relatively small due to the limitation of the hardware. Therefore, the multi-scale method is one of the methods to solve the above dilemma. 
	
	For numerical methods of structural fracture modeling, there are many effective and wild-used numerical  methods at both micro-scale and macro-scale mentioned above, but the coupling between macroscopic continuum model and micro atomic discrete model is the obstacle to develop of multi-scale methods. To solve the problem, many methods have been proposed. Mullins proposed the FEAt method, which corresponds nodes to atoms one by one. This method is simple and easy to use, but it does not consider the influence of long-range forces between atoms, so the force on the boundary atom is unbalanced  \cite{mullins1982simulation}. Then Kohlhoff  improved the method and obtained good results in simple structural problems by use the the displacement connection at the coupling domain instead of  the force connection. \cite{kohlhoff1991crack}. Although   FEAt  method is simple and easy to operate, but it usually requires a large amount of computation. Therefore,  Tadmor  et al.  proposed   QCC (Quasicontiunuum) method to reduce the computational cost by building a  typical atomic system and connecting it to the element of the finite element method (FEM) \cite{tadmor1996quasicontinuum}. The method has been widely used in the field of crack propagation  \cite{tadmor2003peierls} and has obtained good results in two-dimensional plane problems. The MD-FEM multi-scale is a method that combines molecular dynamics (MD) with FEM. Usually, FEM is used to calculate the displacement and stress results in the region far away from the defect while MD is carried out in the area near the defect. This method takes into account both macroscopic and microscopic aspects, which can easily obtain macroscopic parameters such as stress intensity factor (SIF) and analyze the evolution mechanism of microscopic defects. Rafii - Tabar et al. proposed the  MD-FEM method to simulate the two dimensional   Ag   plate crack propagation and has analyzed the propagation of dynamic crack in multiple scales \cite {rafii1998multi, rafii1998multi2}. In addition to the above multi-scale methods, another multi-scale method, named MD-CFEM, is widely used to analyze the behavior of cracks, which combined with the Cohesive Finite Element Method (CFEM) and MD. The MD-CFEM was first proposed by  Yamakov  et al. to investigate the propagation and evolution mechanism of cracks at different scales  \cite{yamakov2006Molecular} . They used MD to simulate the crack propagation of the two-crystal  Al  atom model and extracted the  T-S  curve from it. Then, the parameters of  T-S  curve are transmitted into the  CFEM to simulate the propagation of macro cracks. After that,  MD-CFEM has also been used to simulate the inter-facial fracture behavior of aluminum-based ink coating by Jiang et al. \cite{jiang2018molecular}. Elkhateeb et al. simulated the inter-facial fracture behavior of  $\rm Ti_6Al_4V/TiC$  at different temperatures  \cite{elkhateeb2018molecular} and Sazgar et al. analyzed the interface failure behavior of  $\rm Al/Al_2O_3$  composites  \cite{sazgar2017development} by MD-CFEM. Although the above studies have transmitted microscopic characteristics to macroscopic simulation, however, the research is still focus on the macro crack analysis by CFEM while MD is only used to calculate the parameters of  T-S  curve, and the micro-crack initiation can not be realized.
	
	In this study, a multi-scale method, named multi-grid sampling multi-scale (MGSMS) method, is proposed for crack simulation, which is combined with FEM, XFEM and MD. The process of crack propagation is analyzed comprehensively from macroscopic and microscopic scales, where both the microscopic initiation and macroscopic propagation of crack are studied by MGSMS method.
	
	The remainder of this study is organized as follows. The details of the proposed multi-scale method and model for the crack are introduced in section 2 and 3. The results and discussions of multi-scale simulations can be found in section 4. Finally, some conclusions are summarized in section 5.
	
	\section{multi-scale model} \label{section:ms-model}
\begin{figure}[htb]
	\vspace{\baselineskip}
	\centering
	\includegraphics[width=0.25\textheight]{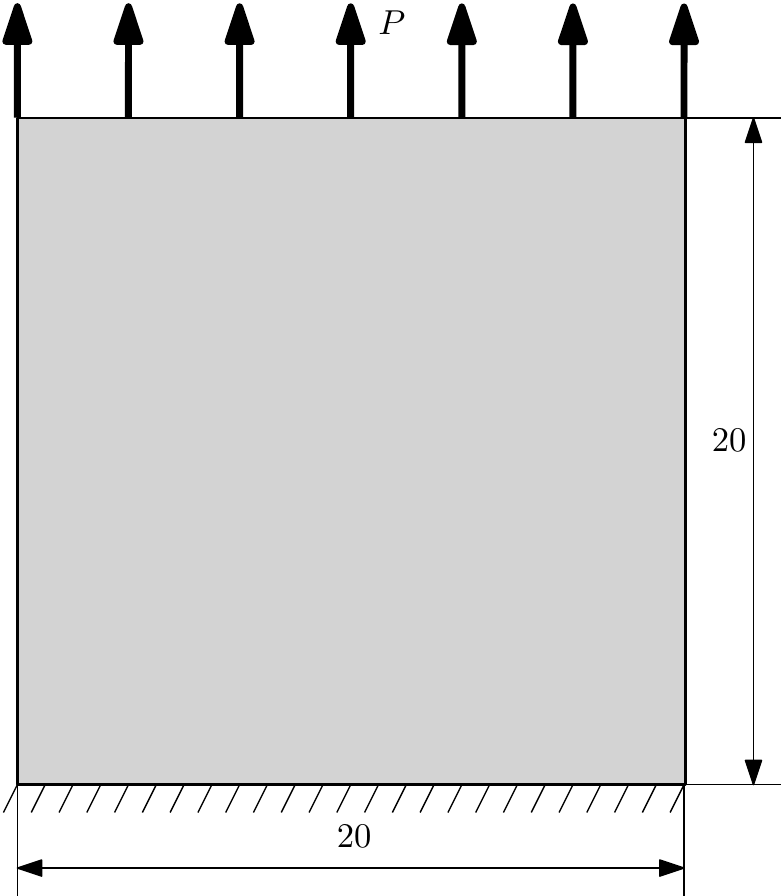}
	\caption{The 2D model of the thin plate  } \label{fig:ms-plate-model}
	\vspace{\baselineskip}
\end{figure}

\begin{figure}[htb]
	\vspace{\baselineskip}
	\centering
	\includegraphics[width=0.6\textheight]{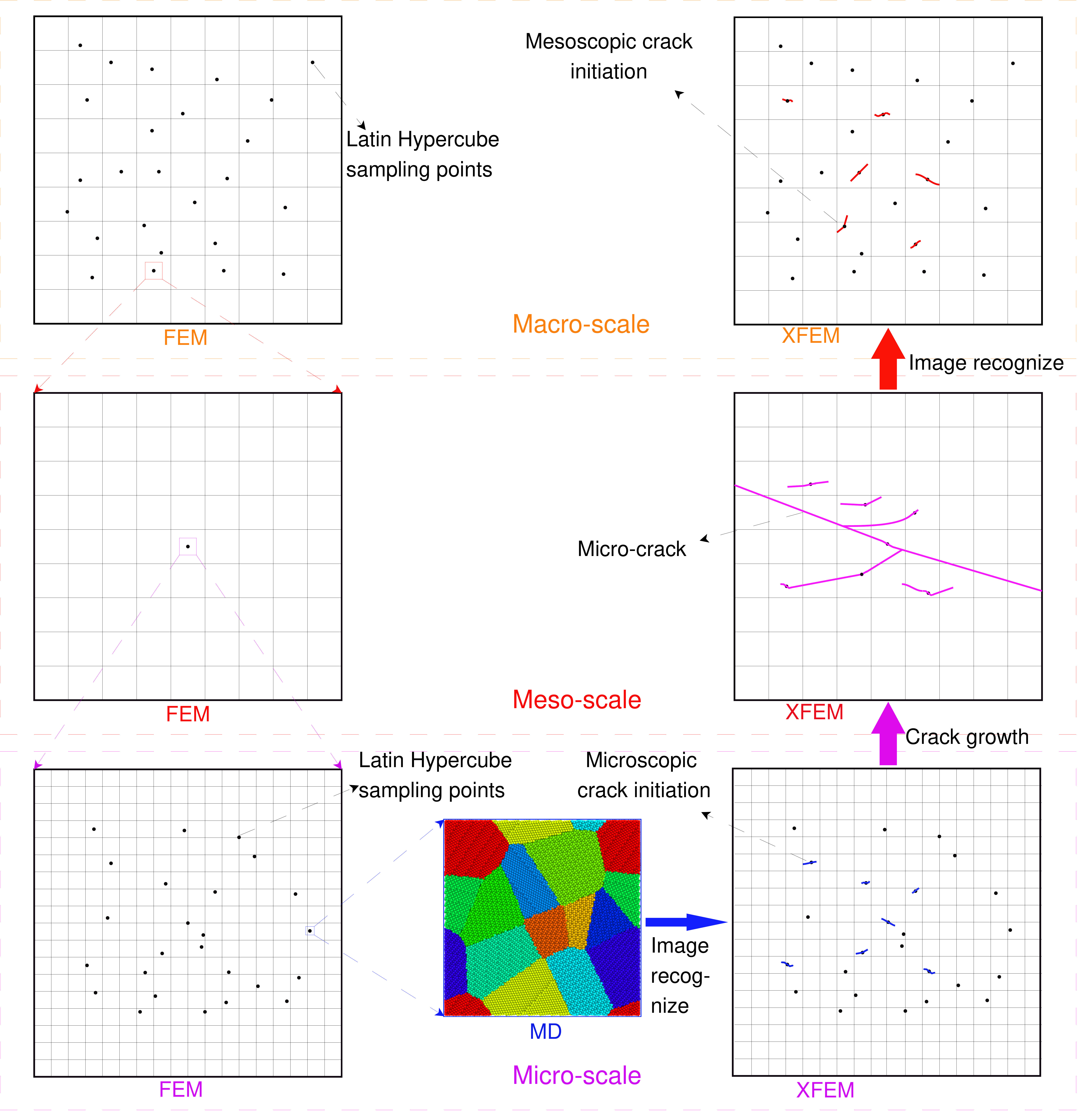}
	\caption{The schematic diagram of the multi-scale} \label{fig:multiscale-model}
	\vspace{\baselineskip}
\end{figure}

	As shown in  Fig. \ref{fig:ms-plate-model} , it is a 2D plate with the size of $20mm \times 20mm$. The lower surface of the plate is fixed while the upper surface is loaded by uniform load  $P $ , where $P=500N/mm $. The material of the plate is low carbon steel, and the parameters are shown in the Tab.  \ref{tab:material}. Assuming that there are no preset cracks on the plate, a novel multi-scale method is proposed to calculate the location of crack initiation and the path of crack propagation after initiating in this study which is combined with FEM, XFEM and MD method. As shown in Fig. \ref{fig:multiscale-model}, FEM and MD are coupled to calculate the crack initiation location by establishing models at macroscopic, mesoscopic and microscopic scale. Then the crack propagation path after crack initiating can be simulated by XFEM. It is obviously that Fig.  \ref{fig:multiscale-model}  can be divided into upper, middle and lower layers, which mean macroscopic, mesoscopic and microscopic models respectively. Moreover, Fig. \ref{fig:multiscale-model}  can be divided into left and right columns. On the left side, the macroscopic displacement boundary conditions are transmitted to the microscopic atomic model layer by layer. On the right side, the crack characteristics obtained from the microscopic MD model are fed back to the macro-scale model step by step. In order to reduce the accuracy error caused by the multi-scale calculation spanning  6 orders of magnitude from the millimeter to nanometer, and to reduce the computational cost caused by the extreme refinement of the mesh, a multi-grid FEM is adopted to extract the displacement boundary conditions in this study. Furthermore, in order to improve the efficiency of the algorithm, only sampling points generated by a Latin hypercube sampling method are used to study the behavior of microscopic crack.
	
	\begin{table*}[htb]
		\small
		\vspace{\baselineskip}
		\centering
		\caption {The material parameters of low carbon steel}
			\label{tab:material}
			\begin{tabular}{c c c c c c}
				\toprule
				$E(Gpa)$ & $\mu$ & $G(GPa)$ & $\sigma_c(MPa) $ &$ m $ &$ C $\\ \midrule
				206&0.3&80&235&2.75&$1.43 \times 10^{-11} $\\ \bottomrule
			\end{tabular}
			\vspace{\baselineskip}
		\end{table*}
		
As shown in Fig. \ref{fig:multiscale-model}, the model size is  $20 mm \times 20 mm $ and the mesh size is  $0.2 mm \times 0.2 mm $ in the macroscopic region. In this step, 50 sample points are generated by Latin hypercube sampling method and the  mesoscopic displacement boundary conditions at each sampling points should be extracted. In the mesoscopic region, the corresponding mesoscopic models need to be constructed according to the displacement boundary conditions of 50 sampling points. The model size in this region is  $200  \mu m \times 200  \mu m $ , and the mesh size is  $2 \mu m \times 2 \mu m $. Then FEM is used to obtain the displacement boundary conditions at smaller scales and transmit them to the microscopic model. In order to transmit the displacement boundary conditions to the atomic discrete model, the model is further refined in the microscopic scale region, where the model size is  $2 \mu m \times 2 \mu m $ and the mesh size is  $20 nm \times 20 nm $. Then the molecular dynamics model is embedded in it. The size of MD model is  $200 \AA \times 200 \AA \times 20 \AA$ . In order to reduce computational cost of MD, Latin hypercube sampling method was used to generate  50  sample points, and then MD analysis was carried out at sample points. After that, the calculation of crack initiation has been completed, and then the propagation behavior after crack initiation needs to be studied. After obtaining the MD simulation result of all the sampling points, it is necessary to extract the position and shape of micro crack initiation from the MD result files. An image recognition method is used in this study. Then, XFEM is used to calculated the path of micro crack propagation according the microscopic crack initiation and then the mesoscopic crack initiation can be obtained from the results of XFEM. Finally, XFEM is used to calculated the path of macro crack propagation according the mesoscopic crack initiation. Thus, the calculation from micro crack initiation to macro crack propagation is completed.
	
\section{multi-scale method}
This study proposed a multi-grid sampling based multi-scale (MGSMS) method for crack simulation. First, the multi-grid based FEM is used to extract the essential boundary conditions from the macro finite element (FE) model to micro atomic  model layer  by layer. Then the multi-grid based XFEM is used to calculate the crack path after initiation and transmit the microscopic crack initiation to the macro model. The process of MGSMS mthod is shown in Fig. \ref{fig:multiscale-model}, and the following is the detail of MSDMS method.
	
\subsection {Coupling boundary condition of FEM and MD} \label {section:fem-md}
	In multi-scale methods, there are two main ways to establish the relationship between the computational methods of macroscopic or mesoscopic continuum and microscopic discrete particles: displacement and force. In this study, the way of  displacement is used to connect the nodes of FE and the particles of MD in the coupling region, which greatly reduces the error caused by the transmission of boundary conditions in the coupling region and avoids the phenomenon of "ghost force" in the coupling region.
	
\begin{figure}[htb]
	\vspace{\baselineskip}
	\centering
	\includegraphics[width=0.6\textheight]{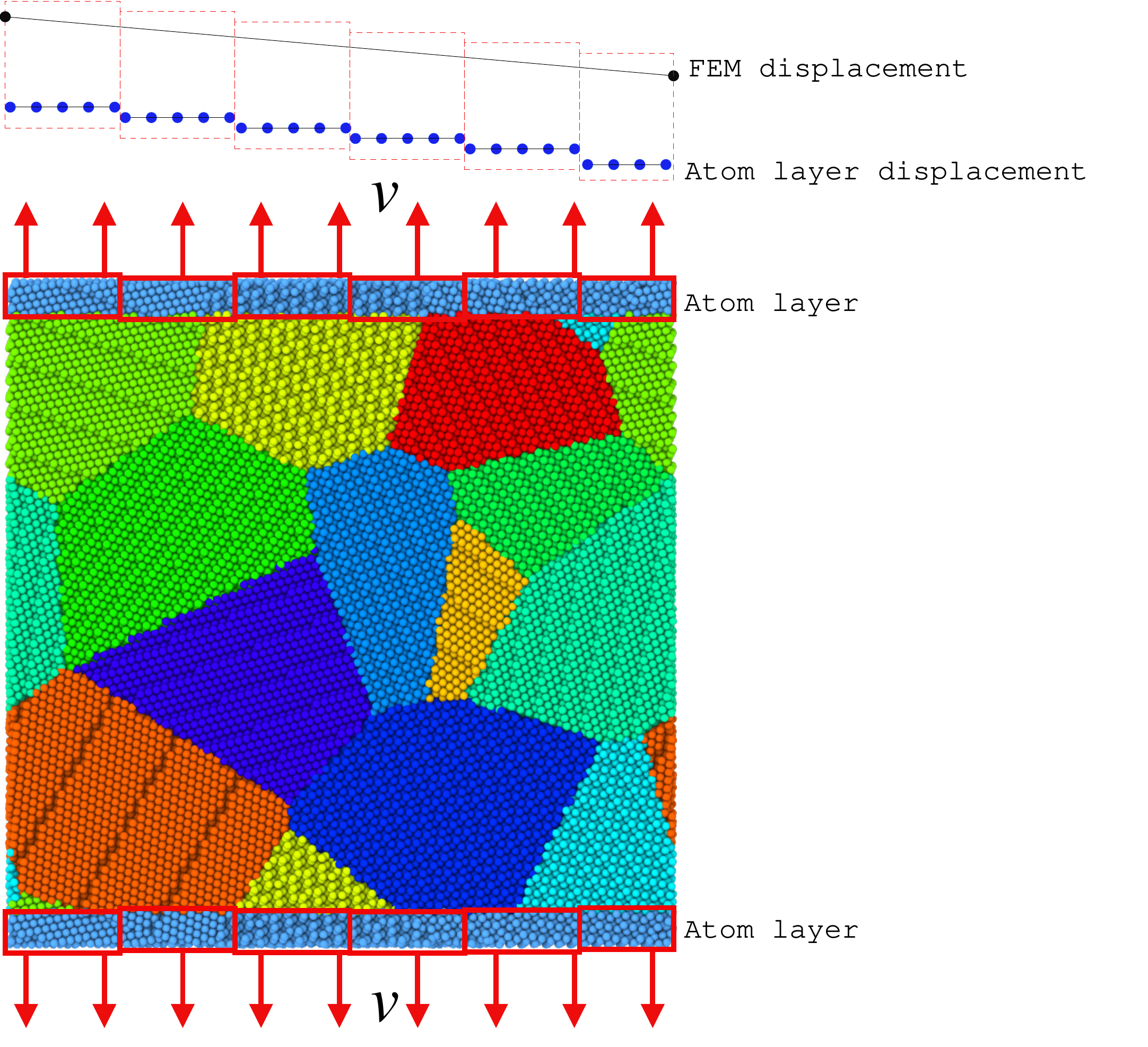}
	\caption{The  schematic diagram of the displacement coupling mode at the junction of FEM and MD} \label{fig:multiscale-load}
	\vspace{\baselineskip}
\end{figure}

	Figure  \ref{fig:multiscale-load}  is a schematic diagram of the displacement coupling mode at the junction of FEM and MD. Firstly, the FEM displacement of the nodes in the coupling region should be extracted from the whole FEM  displacement results. Then the extracted FEM displacement should be divided into several segments and applied to the boundary atomic layer of the MD model. Obviously, the more segments are used, the more accurate the boundary conditions of the MD model will be. However, some scholars have studied this problem in detail   \cite{wang2018ms}. The result shows that the linear loading of  4  and  6  segments have only little effect on the computational results of MD, and the further subdivision has little effect on the results. Therefore, the  6  segments linear loading method is adopted in this study to apply the FEM boundary condition to the MD model. As shown in  Fig.\ref{fig:multiscale-load} , fixed atomic layers with thickness of  $15 \AA $  are set on the upper and lower surfaces of the atomic model for loading the upper and lower boundary conditions. The upper and lower atom layers will be divided into  6  segments. Then, the velocity components in X and Y axes should be calculated according to the FE displacement results at the boundary and the time step of MD simulation. Finally, the calculated velocity components are assigned to each atomic layer to applied the displacement boundary conditions.
	
	\subsection {Molecular dynamics simulation in discrete atomic region} \label {section:ms-md-cal}
	
\begin{figure}[htb]
	\vspace{\baselineskip}
	\centering
	\subfigure{\includegraphics[width=2in]{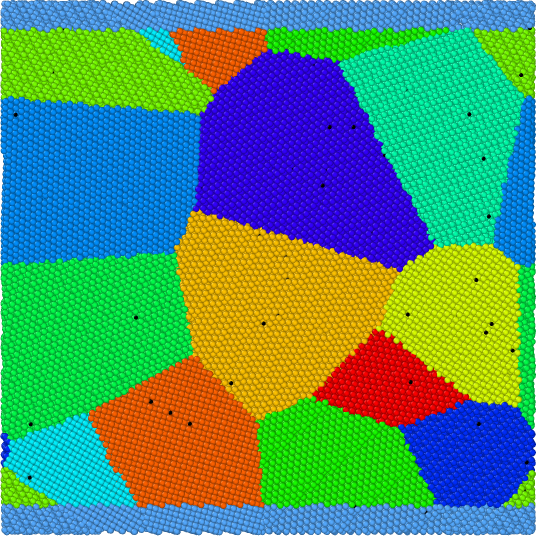}}
	\subfigure{\includegraphics[width=2in]{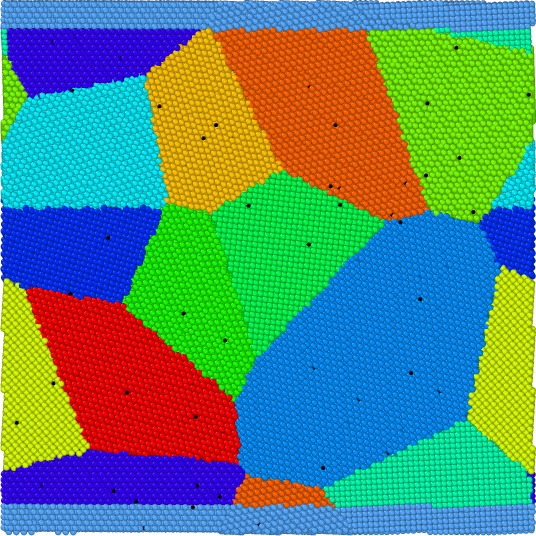}}
	\subfigure{\includegraphics[width=2in]{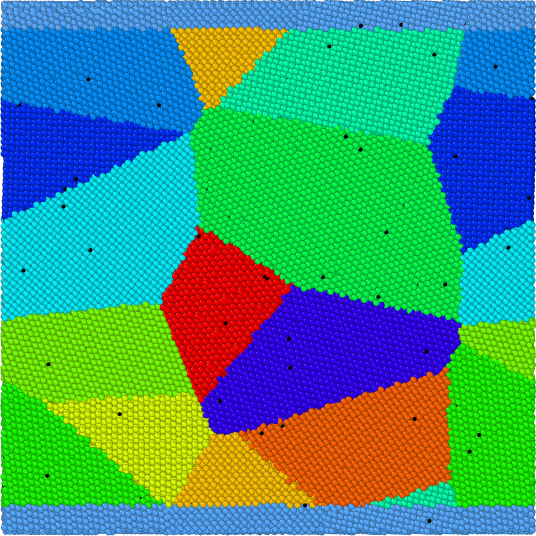}} 
	\caption{The poly-crystal atomic models generated by Atomsk} \label{fig:ms-md-models}
	\vspace{\baselineskip}
\end{figure}

	As shown in Fig.  \ref{fig:multiscale-model} , it can be found that the MD method is used to simulate the crack initiation at the micro scale in the process of multi-scale computation. In this study, atomic models of microscopic crystals are constructed according to the  Voronoi diagram which is generated by algorithms randomly to reflect the complexity and diversity of micro-structure. The examples of atomic models are shown in  Fig. \ref{fig:ms-md-models}. It can be found that the microscopic crystal structure of each sampling point is with different grain arrangement and orientation. This study takes steel as an example to study the process of micro-crack initiation by constructing the micro-crystalline structure of steel. In this study, a Fe atoms polygonal crystal model is constructed with  $0.2\% $  vacancy and  $0.2\% $  carbon atoms in which each single grain is arranged with  BCC  crystal structure. The lattice constant is  $a=2.856 \AA $ . The crystal orientation of each grain is randomly generated by the algorithm. The size of the model is  $200 \AA \times 200 \AA \times 20 \AA $ , and it contains about   68000   atoms inside. An open source software  named Atomsk\cite{hirel2015atomsk}  is used to generate atomic models. Since there are hundreds of  atomic models need to be generated, so all the models can be generated by a script. The modeling script is as follows:
	
\begin{lstlisting}[language=c++]
	#### Polycrystal model with defects by Atomsk  ####
	atomsk --create bcc 2.855 Fe fe.xsf
	# Polycrystal
	atomsk --polycrystal fe.xsf fe_poly.txt fe_poly.cfg -wrap
	# Vacancies
	atomsk fe_poly.cfg -select random 0.2% any -rmatoms select fe_poly_v.cfg
	# Substitude atoms
	atomsk fe_poly_v.cfg -select random 0.1% Fe -substitute Fe C fe_poly_v_c.cfg
	# Interstitial atoms
	atomsk fe_poly_v_c.cfg  -add-atom C random 68 fe_poly_v_c.cfg .cfg lmp
\end{lstlisting}
	
	Generally, atomic models built directly are often in a non-equilibrium state due to unstable surface energy. Therefore, it is necessary to relax the model at room temperature (300K) first to reach the initial equilibrium state. After completing the relaxation process,  the next step is to apply the boundary conditions to simulate the internal crack initiation. Fixed atomic layers with thickness of  $15\AA$  are set on the upper and lower surfaces of the model respectively for loading the upper and lower boundary conditions. How to apply the boundary conditions to the atomic model can be found in section \ref{section:fem-md}. Moreover, the time-step of MD simulation is set as  $0.001PS $. The simulation ensemble is set as NVT and the temperature of  NVT  ensemble is set as  300K. The potential function is selected as  EAM  potential \cite{Hepburn2008metallic}. The boundary conditions of the model are set as flexible boundary conditions in X and Y  axes and periodic boundary conditions in Z axis. In this study,  LAMMPS \cite{plimpton1995fast} is used to execute MD simulations. All the process of MD simulations including initialize, set input and output and running control are operated by a LAMMPS  script file.
	
	\subsection{Image recognition based crack extracting method} \label{section:image-process}
	
	In the computational process of LAMMPS, it is difficult to output the real-time position of the crack directly due to the lack of crack recognition mechanism. Moreover, it will take in a long time for data processing when extracting the position of crack from the result file. Because the size of LAMMPS result file is usually very large (several gigabytes for system with 68000 atoms). Therefore, an image recognition based crack extracting method is proposed in this study. Image recognition technology can  extract the required feature from the picture intuitively. It is an effective method for crack recognition in this study. Moreover, since the result files of several gigabytes (GB) or several decades  GB are converted into images with only several hundred kilobytes (KB) , so the image recognition based crack extracting method can not only save a large amount of storage space is saved, but also significantly reduce the crack extracting time.
	
\begin{figure}[htb]
	\vspace{\baselineskip}
	\centering
	\includegraphics[width=0.6\textheight]{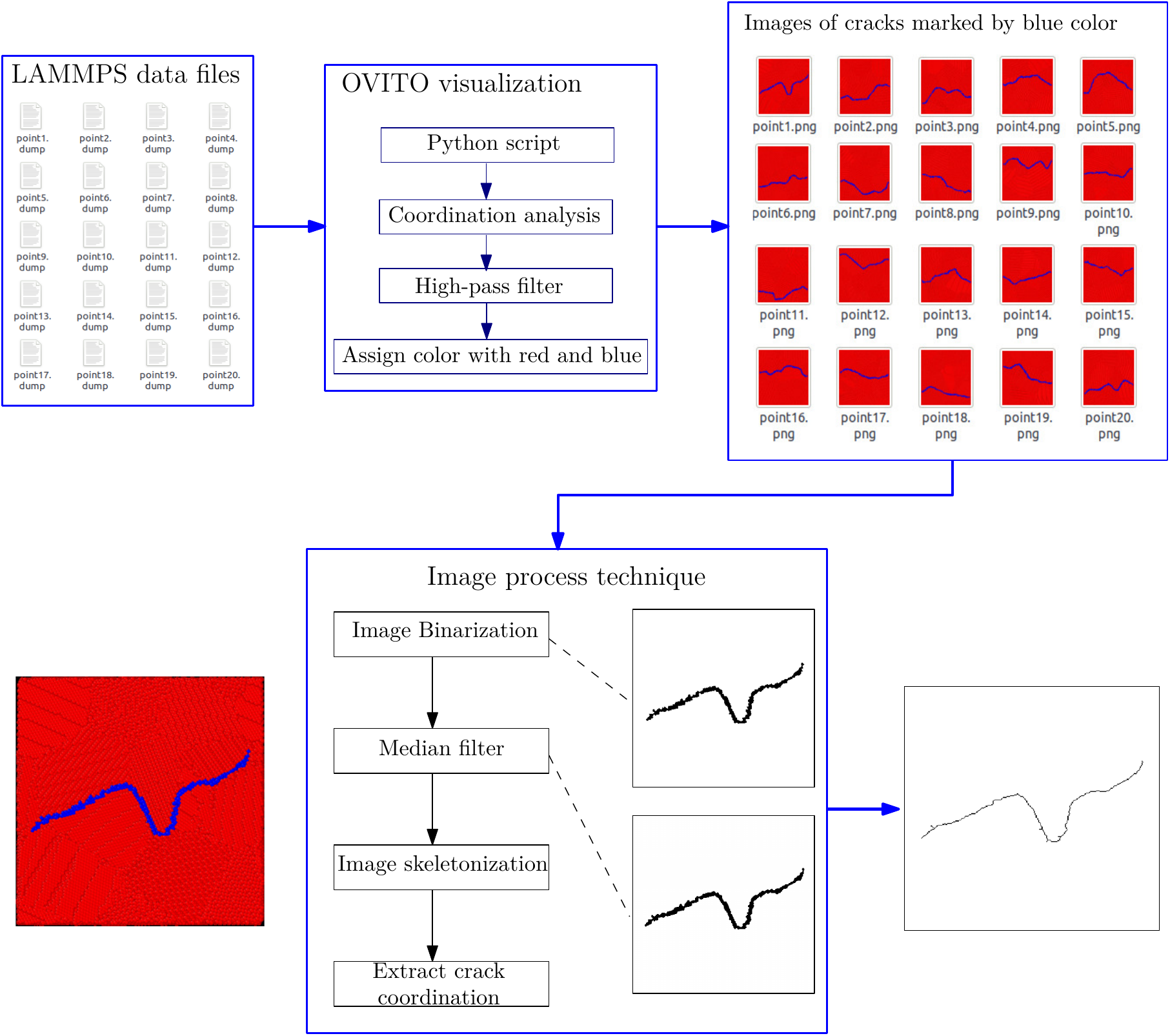}
	\caption{The flowchart of the image recognition based crack extracting method} \label{fig:ms-image-process1}
	\vspace{\baselineskip}
\end{figure}

	Figure \ref{fig:ms-image-process1}  is the flowchart of the image recognition based crack extracting method. It can be found that the first step of this method is to use the post-processing software  (OVITO  \cite{ovito}) to process the result files of LAMMPS  and output the corresponding result picture. OVITO  is a software which is used for visualization and analysis of atomic and particle simulation data. Its built-in crystal structure analysis, dislocation, defect statistics and other functions provide strong support for the analysis of MD results. Moreover, OVITO  also supports command-line operation, which automatically analyzes the result file and generates the desired result image by calling  Python  script. In this section, the coordination analysis function in OVITO  is used to identify the crack location in polygonal crystals. This function is mainly used to calculate the number of adjacent particles of each particle within a given cutoff radius around its position, known as the coordination number. Obviously, the coordination number will be relatively small at the model boundary or fracture surface. Then a high-pass filter can be used to separate the model boundary or fracture surface region from other regions. After that, because the model boundary region is easy to exclude, so the fracture surface region can be isolated. Then, atoms in the isolated fracture surface region will be assigned with blue color while the others are assigned with red color. At this point, the pre-processing work for crack extracting has been ready, and all the above operations are executed by calling the written  Python  script in  OVITO , which can generate the result pictures in batch. The  Python  script used in this study is as follows:
	
	\begin{lstlisting}[language=python]
		# Python script for OVITO 3.0.0 to generate images from LAMMPS output files
		# Zhenxing Cheng, 2021
		# Email: zxcheng@hnu.edu.cn
		
		from ovito.io import *
		from ovito.modifiers import *
		from ovito.pipeline import *
		from ovito.vis import *
		# Data import:
		pipeline = import_file('./output/fe_poly_mmp_t0.dump', multiple_frames = True)
		pipeline.add_to_scene()
		# Visual appearance setup:
		pipeline.source.data.cell.vis.render_cell = False
		pipeline.source.data.particles.vis.radius = 2.0
		# Coordination analysis:
		pipeline.modifiers.append(CoordinationAnalysisModifier(
		cutoff = 7.6, 
		number_of_bins = 100))
		# Freeze property:
		pipeline.modifiers.append(FreezePropertyModifier(
		source_property = 'Position', 
		destination_property = 'Position'))
		# Assign color:
		pipeline.modifiers.append(AssignColorModifier(color = [1.0, 0.0, 0.0]))
		# Expression selection:
		pipeline.modifiers.append(ExpressionSelectionModifier(expression = 'Position.X
		>10 && Position.X<190 && Position.Y>10 && Position.Y<190 && Coordination<115'))
		# Assign color:
		pipeline.modifiers.append(AssignColorModifier(color = [0.0, 0.0, 1.0]))
		# Viewport setup:
		vp = Viewport(
		type = Viewport.Type.Top)
		vp.zoom_all()
		# Rendering:
		vp.render_image(filename='crack.png', size=(330, 330), alpha=True, frame=100,
		 renderer=None)
	\end{lstlisting}

	As mentioned above, after generating the result images from the MD result files by OVITO, it is necessary to use specific image processing techniques to recognize the position of cracks from images. The first step is to transfer result images into binary images with only  0  and 1. Then the median filter is used to eliminate the noise in the image to ensure the accuracy of results. Subsequently, image skeletonization processing is used to extract the skeleton line of the crack. Then, the relative coordinate value  $(x_r, y_r) $of the crack is extracted from the skeleton line of the crack, where the coordinate value is calculated based on the pixel points of the image. The physical coordinate value of the crack  $(x,y)$ can be obtained by Eq.  \eqref{eq:xcoord}  and Eq. \eqref{eq:ycoord}, where $l_x $ and $l_y $  mean the size of atomic model in the X and Y axes respectively. Symbols  $n_x $ and $n_y $ are the number of pixel of result image in X and Y axes respectively. Symbols $x_0 $  and  $y_0 $  are the actual physical coordinates of the reference point (the vertex in the lower left corner).
	
	\begin{gather}
		x = x_r (\frac{l_x}{n_x})+x_0 \label{eq:xcoord} \\
		y = y_r (\frac{l_y}{n_y})+y_0 \label{eq:ycoord}
	\end{gather}

\subsection {Calculation of multiple extended finite element method}

As shown in Fig.  \ref{fig:multiscale-model} , it can be found that after the micro-crack initiations are obtained, the multi-grid based XFEM should be used to simulate the crack propagation after crack initiating and then feedback the micro-crack initiations to the macroscopic model step by step.  Firstly, the MD results of  50 Latin hypercube sampling points are analyzed at the microscopic scale, and the points which contain crack initiation will be screened out. Then, the crack coordinates of each sample point should be extracted by the image recognition based crack extracting method. Then, a micro propagation finite element model was constructed according to coordinates of crack initiations, and XFEM is used to simulate the multiple cracks propagation under the micro model by applying the corresponding displacement boundary conditions. In order to simulate the phenomenon of asynchronous propagation of multiple cracks, the theory of maximum energy release rate of crack propagation is used in this study  \cite{nuismer1975energy}. This means that only part of the micro-crack initiations will grow further and lead to the complete fracture of the model eventually. After all the microscopic cracks stop growing, the main crack is selected as the initiation of mesoscopic crack according to the final crack propagation path and transmitted it to the macro model. Subsequently, the crack coordinates of the mesoscopic crack initiations of each sampling point can be determined according to the mesoscopic crack simulation results of  50 Latin hypercube sampling points. Finally, the macroscopic extended finite element model is constructed based on the results of mesoscopic crack propagation, and the XFEM should be used to calculate the final macroscopic crack propagation path.

\section{results and analysis}
As shown in  Fig. \ref{fig:ms-plate-model} , a plate with uniformly distributed loading is taken as an example to analyze the internal crack initiation and propagation by MGSMS method. The model size and material parameters involved in this section have been given above. In this section, only the calculated results are given and necessary analysis is made.

\subsection {Simulation results of multi-grid based finite element method}

\begin{figure}[htb]
	\vspace{\baselineskip}
	\centering
	\subfigure{\includegraphics[width=0.3\textheight]{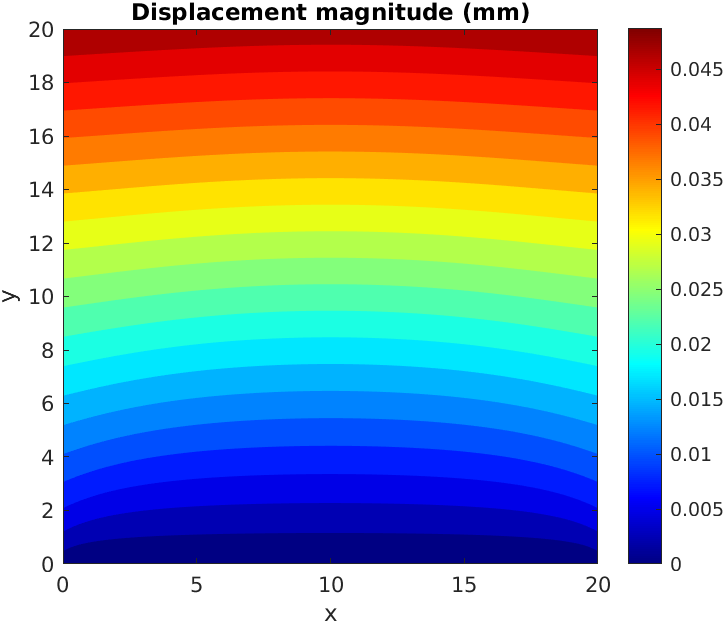}}
	\hspace{0.1in}
	\subfigure{\includegraphics[width=0.3\textheight]{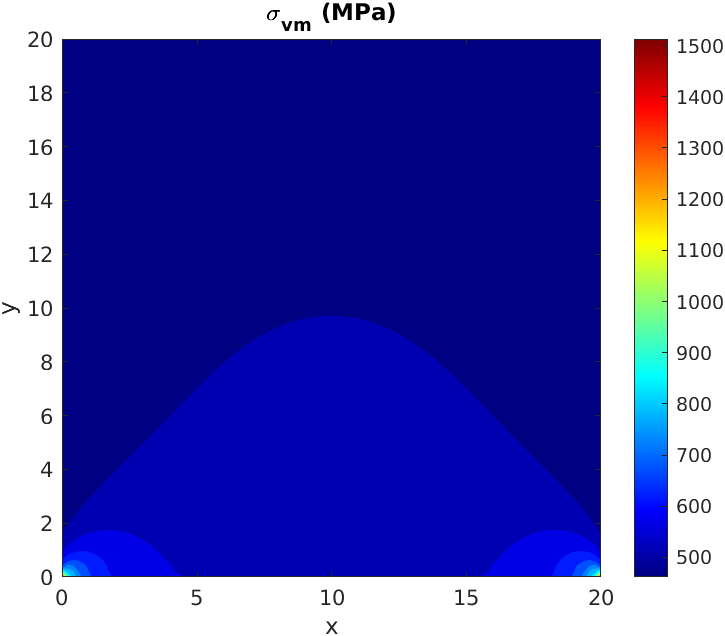}}
	\caption{The analysis result of FEM at macroscopic scale} \label{fig:ms-fem-mac-results}
	\vspace{\baselineskip}
\end{figure}

In this study, the macroscopic deformation is transmitted to the microscopic atomic model through the displacement boundary conditions by the multi-grid FEM. In order to meet the precision requirements of the micro-scale analysis, three-level FEM (macro, mesoscopic and microscopic) is used to amplify the macroscopic displacement results layer by layer. Figure  \ref{fig:ms-fem-mac-results}  is the FEM result at macroscopic scale. The left image is the contour plot of the displacement magnitude while the right is the contour plot of the Von Mises stress. At the macroscopic scale, the model size is  $20 mm \times 20 mm $ , and the finite element discrete mesh is a quadrilateral mesh with the size of  $0.2 mm \times 0.2 mm $, and the number of mesh is  $100 \times 100 $.

\begin{figure}[htb]
	\vspace{\baselineskip}
	\centering
	\includegraphics[width=0.6\textheight]{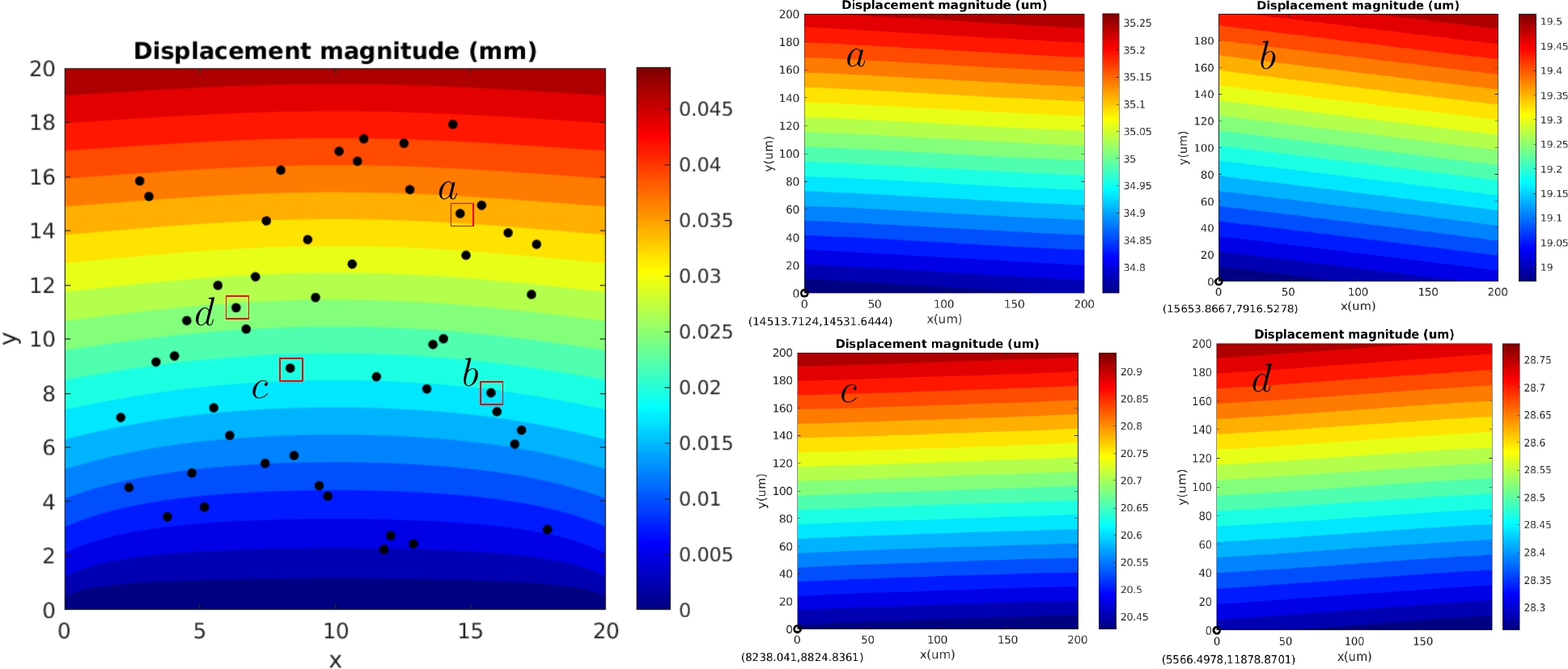}
	\caption{The analysis result of FEM at mesoscopic scale} \label{fig:ms-meso-fem-dis}
	\vspace{\baselineskip}
\end{figure}

In order to obtain the displacement value of the model at mesoscopic scale, it is necessary to refine the mesh further. Then the number of grids can be as high as  $10000\times10000 $  level if all the grids is refined. Howerver, it is difficult to achieve such computational scale on a normal computer. Therefore, the sampling method is adopted to generate sampling points for mesoscopic scale analysis in this study, which greatly reduces the number of grids required for multi-scale analysis. In this study, Latin hypercube sampling method is used to generate  50  sample points, and each sample point should be analyzed at mesoscopic scale by FEM, the analysis results are shown in Fig.  \ref{fig:ms-meso-fem-dis}. The left image is the contour plot of macroscopic displacement, in which the black points are the sampling points. Some examples of the mesoscopic displacement contour at the sampling points are given as shown in the right four images. At this stage, the model size at each sampling point is  $200  \mu m \times 200  \mu m $ and the mesh size is  $2 \mu m \times 2 \mu m $ , which has reached the microscopic scale region.

\begin{figure}[htb]
	\vspace{\baselineskip}
	\centering
	\includegraphics[width=0.6\textheight]{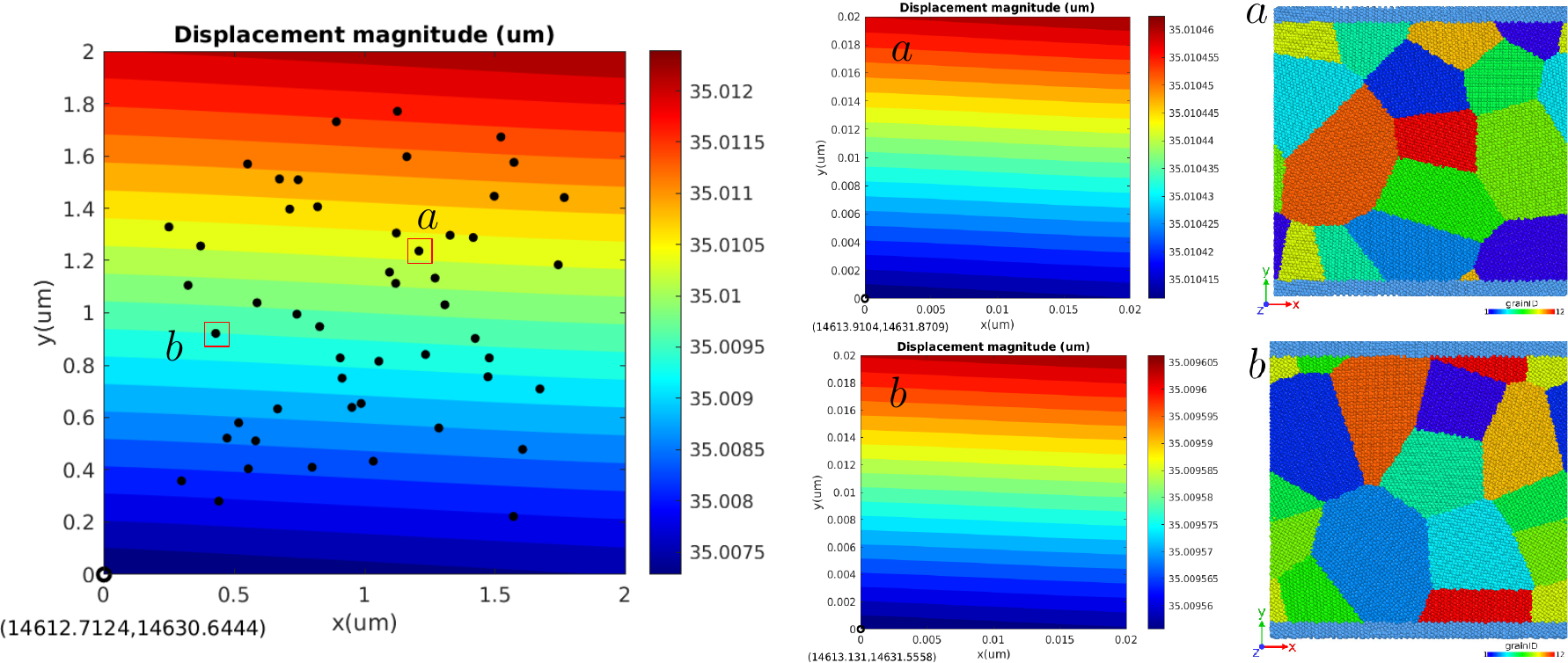}
	\caption{The analysis result of FEM at microscopic scale} \label{fig:ms-micro-fem-dis}
	\vspace{\baselineskip}
\end{figure}

In order to transmit the displacement boundary conditions to the atomic discrete model, the model need to be further refined in the microscopic scale region, where the model size is  $2 \mu m \times 2 \mu m $ , and the FE mesh with $100 \times 100 $  is used for FEM analysis. The size of mesh is $20 nm \times 20 nm $ at this step. Then MD is used to study the crack initiation at the micro scale. Moreover, because the computational process of MD is time consuming, so if the entire region of  model ($2 \mu m \times 2 \mu m $)  was used to build the atomic model for MD simulation, the number of atoms should be up to  $10^9 $  order, which is also difficult to achieve on a normal computer for this computational task. Therefore,  Latin hypercube sampling method is used again to reduce the computational cost in this study by generate  50 sample points. Then MD method is carried out on each sample point. The size of the MD model at each sample point is consistent with the size of mesh ($200 \AA \times 200 \AA \times 20 \AA$). As shown in Fig.  \ref{fig:ms-micro-fem-dis}, the left image is the displacement contour of the micro-scale model, where the black points are the sampling points. Two points  $a $  and  $b $  were selected to explain the computational process. First, extract the displacement contour with the size of a single mesh ($200 \AA \times 200 \AA $) at each sample points. Then construct poly-crystal atomic models with the same size ($200 \AA \times 200 \AA \times 20 \AA$), which is shown in the right two images of Fig. \ref{fig:ms-micro-fem-dis}. The poly-crystal atomic models are constructed according to the Voronoi  diagram generated randomly by Atomsk. Each model contains  12  grains and the detail for how to generate atomic models was introduced in section \ref{section:ms-md-cal}.

In summary, this section described the detail for how to transmit the macroscopic deformation to the microscopic atomic model by displacement boundary conditions. The displacement results are amplified by three-level FEM and simplified by two-level Latin hypercube sampling method. Therefore, 2500 sample points were generated in this study, which means that  2500  different microscopic poly-crystal models were constructed for MD analysis.

\subsection {Simulation results of molecular dynamics }

\begin{figure}[htb]
	\vspace{\baselineskip}
	\centering
	\subfigure[Stress-strain curve]{\includegraphics[width=0.4\textheight]{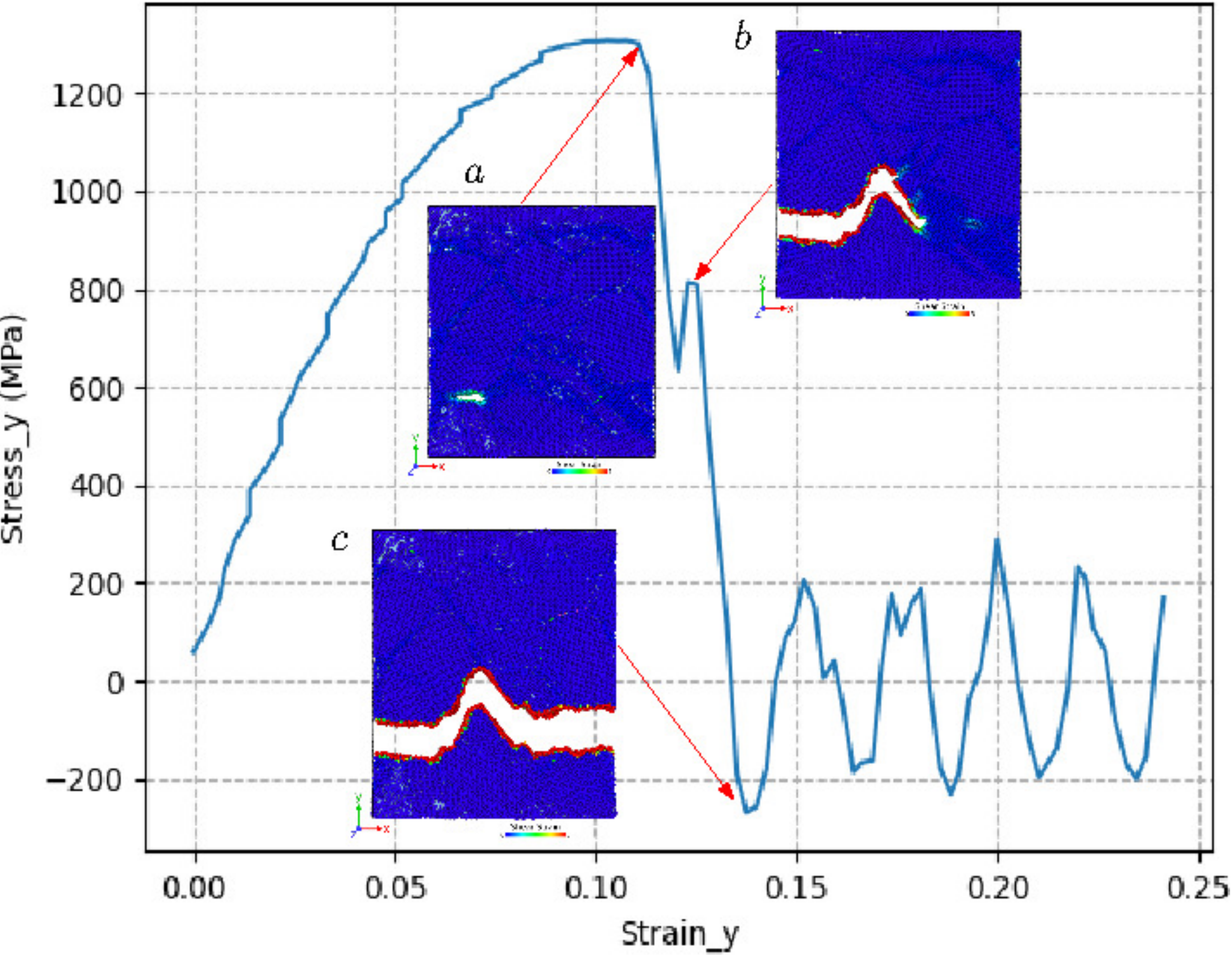}}
	\hspace{0.1in}
	\subfigure[Final fracture diagram]{\includegraphics[width=0.2\textheight]{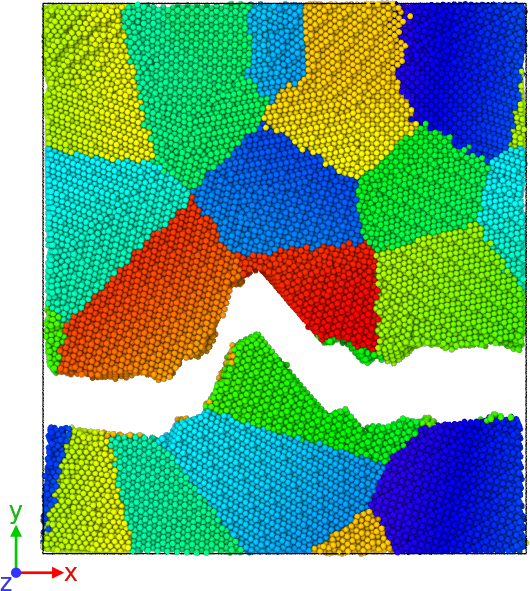}}
	\caption{The stress-strain curve and final fracture diagram of the microscopic poly-crystal at point $ a $} \label{fig:ms-md-point3-point6}
	\vspace{\baselineskip}
\end{figure}

As mentioned above,  2500  poly-crystal models were constructed to analyze the micro-crack initiation locations in this study. The displacement boundary conditions extracted from the microscopic finite element method are applied to the microscopic poly-crystal models according to the method described in  section \ref{section:fem-md} and the corresponding MD simulation is executed. Several sets of results are selected as examples for specific analysis. Figure \ref{fig:ms-md-point3-point6} shows the MD results corresponding to the sampling point  $a $  in the  Fig.\ref{fig:ms-micro-fem-dis}. The left image is the stress-strain curve of the poly-crystal structure when it is loading, where figures  $a $, $b$  and  $c$  are atomic shear strain plots when the strain is  0.11, 0.13  and  0.14 , respectively. It can be found that, the crack initiation begins to appear at the grain boundary when the strain reaches  0.11. After that, the crack will continue to grow along the grain boundary until the fracture, which as shown in Fig.\ref{fig:ms-md-point3-point6}(b). It is obvious that the crystal is  fractured along the grain boundary eventually.

\begin{figure}[htb]
	\vspace{\baselineskip}
	\centering
	\subfigure[Stress-strain curve]{\includegraphics[width=0.4\textheight]{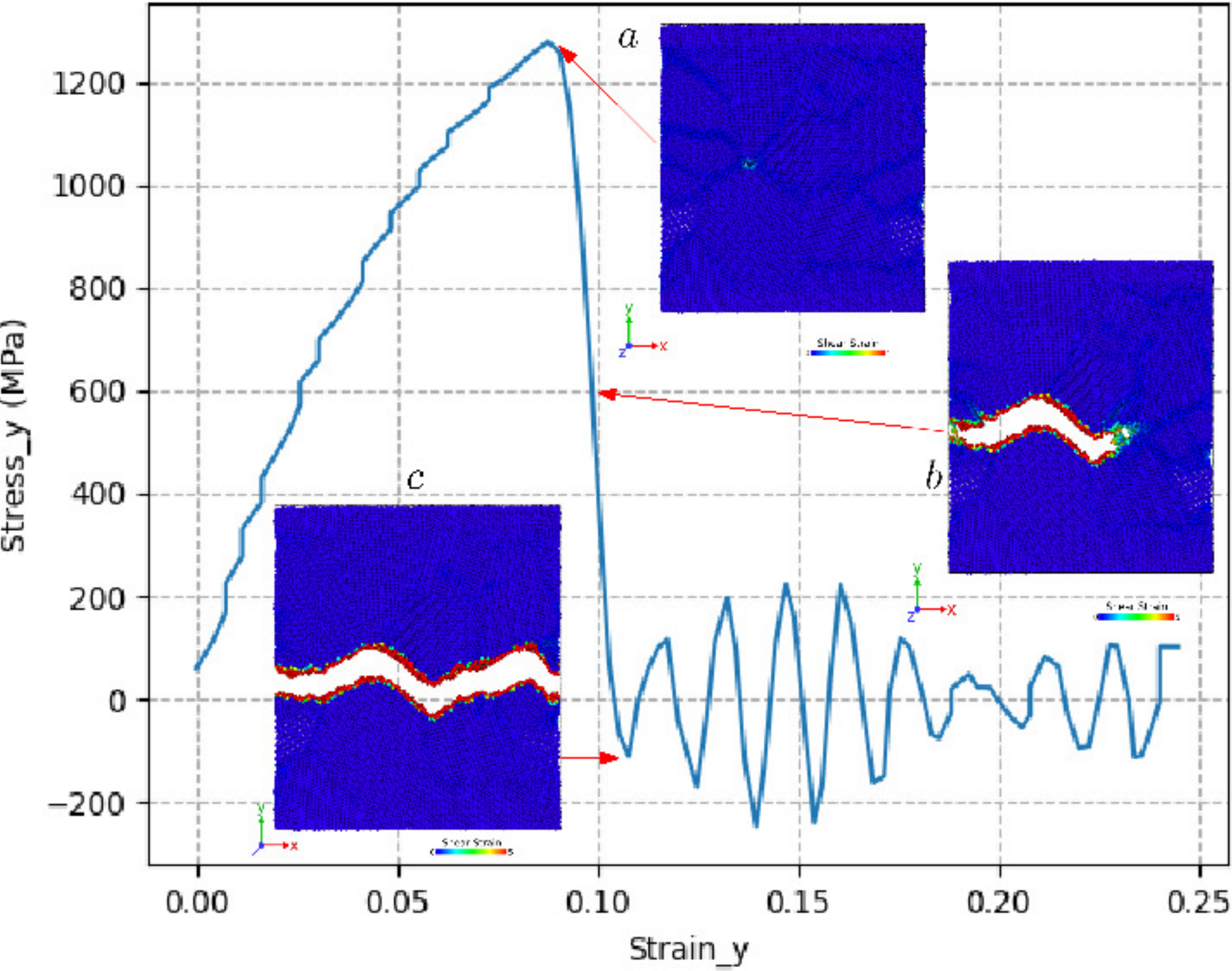}}
	\hspace{0.1in}
	\subfigure[Final fracture diagram]{\includegraphics[width=0.2\textheight]{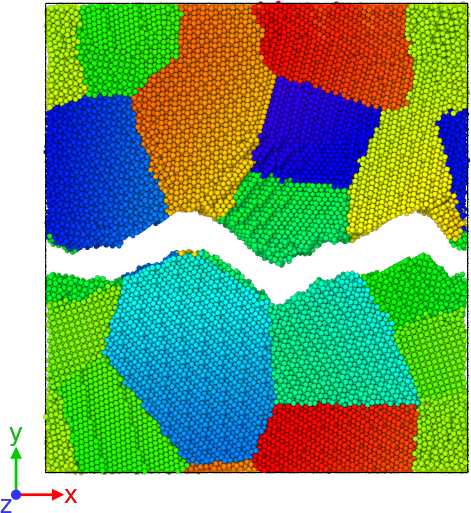}}
	\caption{The stress-strain curve and final fracture diagram of the microscopic poly-crystal at point $ b $} \label{fig:ms-md-point3-point16}
	\vspace{\baselineskip}
\end{figure}

 Another example is the sampling point $b$ in Fig.\ref{fig:ms-micro-fem-dis}. Figure \ref{fig:ms-md-point3-point16}  is the MD result corresponding to the sampling point $b$. The left image also shows the stress-strain curve while the right image shows the final fracture diagram of the crystal.  It can be found that the crack initiation begins to appear in the poly-crystal in this model when the strain reaches  0.09,  and then the crack will grow along the grain boundary rapidly with the continuous loading. By comparing the results between  $a $  and  $b $ , it can be found that the starting time of crack initiation is different for different crystal structures and the crack path is also different, but it usually grow along the grain boundaries.

After completing the MD simulations of  2500  sample points, a series of MD results are obtained. Some of them contains the crack initiations and some of them show no or only partial cracks. Therefore, it is necessary to extract the crack initiations from those MD results files. In this study, the image recognition based crack extracting method introduced in  section \ref{section:image-process} is used to extract the crack initiations from the result files.  

\begin{figure}[htb]
	\centering
	\subfigure{\includegraphics[width=0.19\textheight]{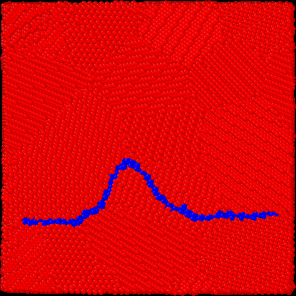}}
	\hspace{0.1in}
	\subfigure{\includegraphics[width=0.19\textheight]{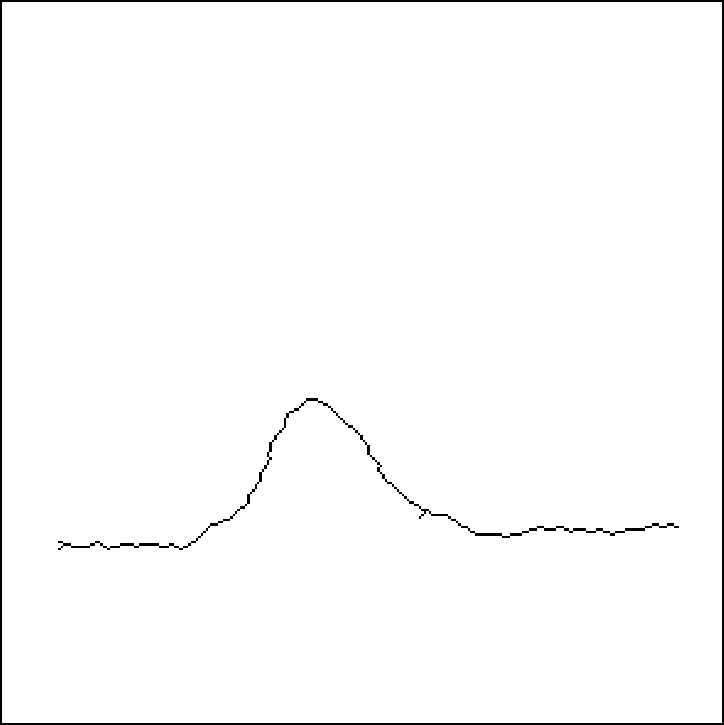}}
	\subfigure{\includegraphics[width=0.2\textheight]{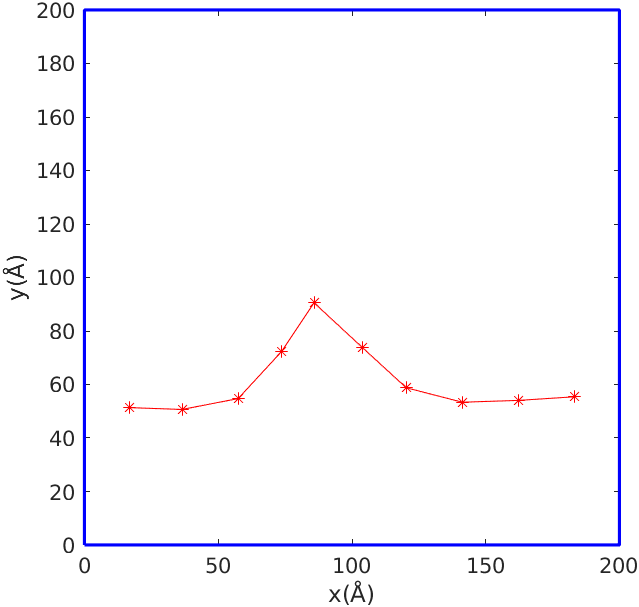}}\\
	\subfigure{\includegraphics[width=0.2\textheight]{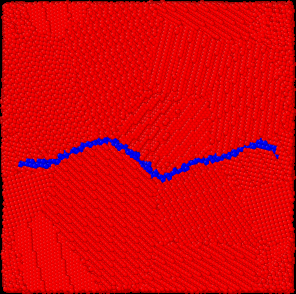}}
	\hspace{0.1in}
	\subfigure{\includegraphics[width=0.19\textheight]{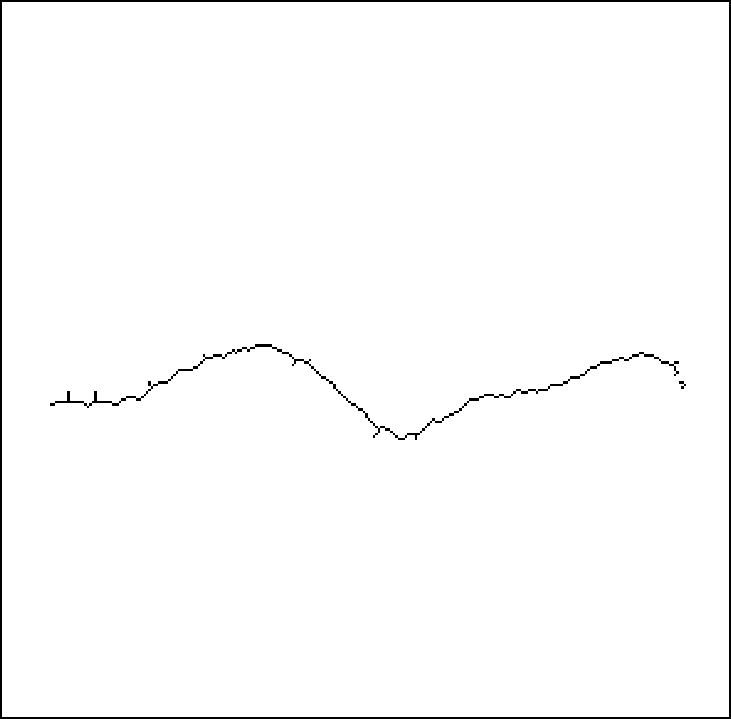}}
	\subfigure{\includegraphics[width=0.19\textheight]{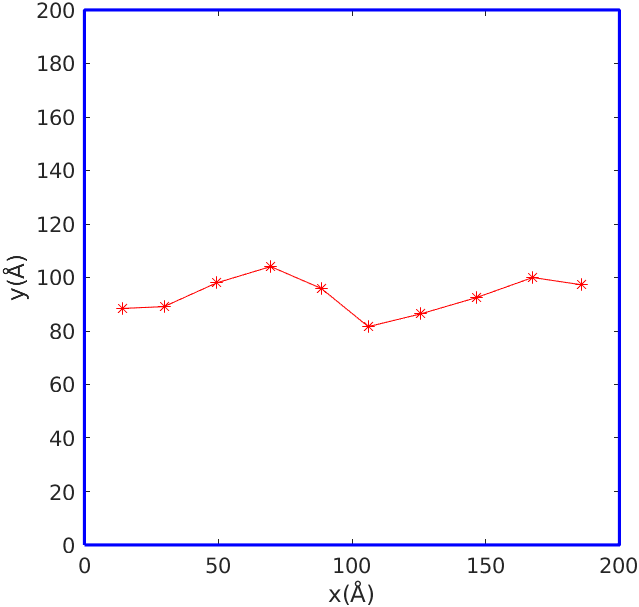}}
	\caption{Extracted crack results at points $ a $ and $ b $ } \label{fig:ms-md-point3-point16-crack-s}
\end{figure}

Figure \ref{fig:ms-md-point3-point16-crack-s}  shows the result extracted by the image recognition based crack extracting method, where $a $  and  $b $  are corresponding to the sampling point  $a $  and $ b $ in the  Fig. \ref{fig:ms-micro-fem-dis} respectively. It can be found that the extracted crack path are conformed to the MD results and the coordinates of extracted cracks are shown in Tab. \ref{tab:crack-coord}.

\begin{table}[htbp]
	\small
	\vspace{\baselineskip}
	\centering
	\caption{The coordinates of cracks extracted by image recognition based crack extracting method}
	\label{tab:crack-coord}
	\begin{tabular}{ccccc}
		\toprule
		\multirow{2}{*}{Number} & \multicolumn{2}{c}{Coordinate of point~$a$~}          & \multicolumn{2}{c}{Coordinate of point~$b$~}          \\ \cline{2-5} 
		& X (\AA)             & Y(\AA)              & X(\AA)              & Y(\AA)              \\ \midrule
		1                   & 16.89 & 51.35 & 14.18 & 88.43 \\
		2                   & 36.48 & 50.67 & 29.72 & 89.11 \\
		3                   & 57.43 & 54.72 & 49.32 & 97.95 \\
		4                   & 73.64 & 72.29 & 69.59 & 104.08 \\
		5                   & 85.81 & 90.54 & 88.51 & 95.92  \\
		6                   & 104.05 & 73.64 & 106.08 & 81.63 \\
		7                   & 120.27 & 58.78 & 125.67 & 86.39 \\
		8                   & 141.21 & 53.37 & 146.62 & 92.517 \\
		9                   & 162.16 & 54.05 & 167.56 & 100.00 \\
		10                  & 183.10 & 55.41 & 185.81 & 97.28 \\ \bottomrule
	\end{tabular}
\end{table}

\subsection {Simulation results of multi-grid based extended finite element method}

\begin{figure}[b]
	\centering
	\subfigure{\includegraphics[width=0.2\textheight]{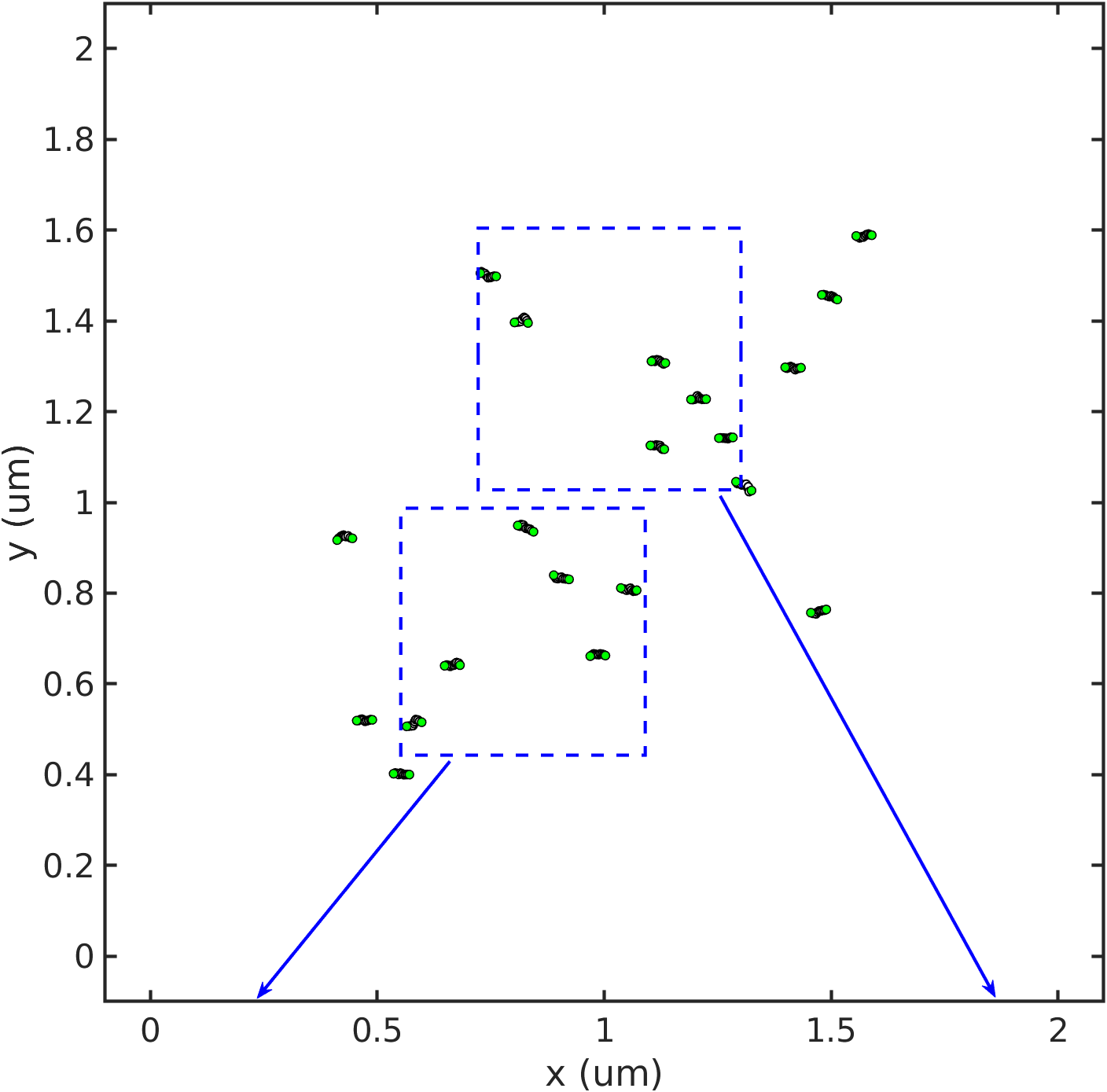}}\\
	\subfigure{\includegraphics[width=0.2\textheight]{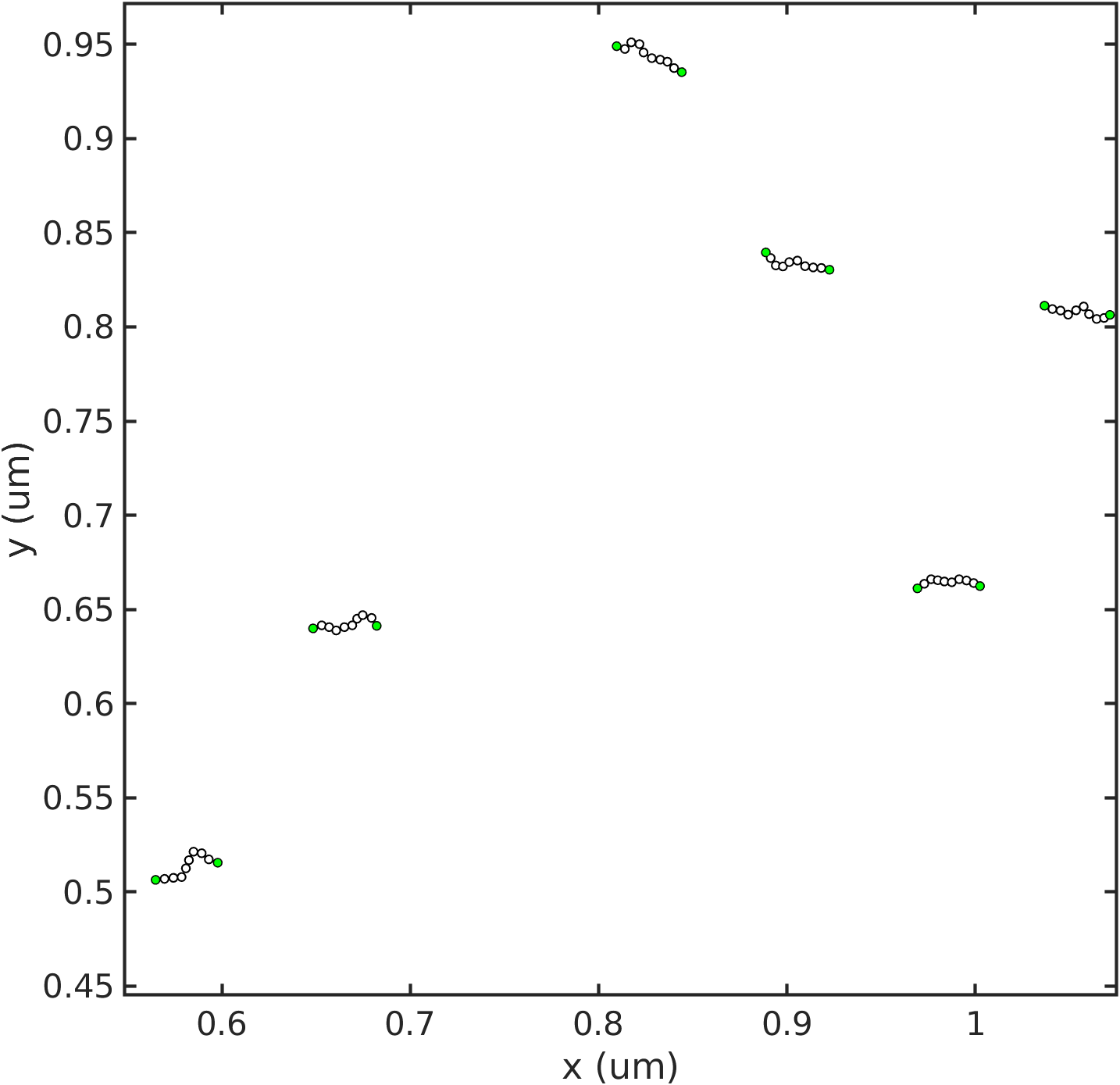}}
	\hspace{0.1in}
	\subfigure{\includegraphics[width=0.2\textheight]{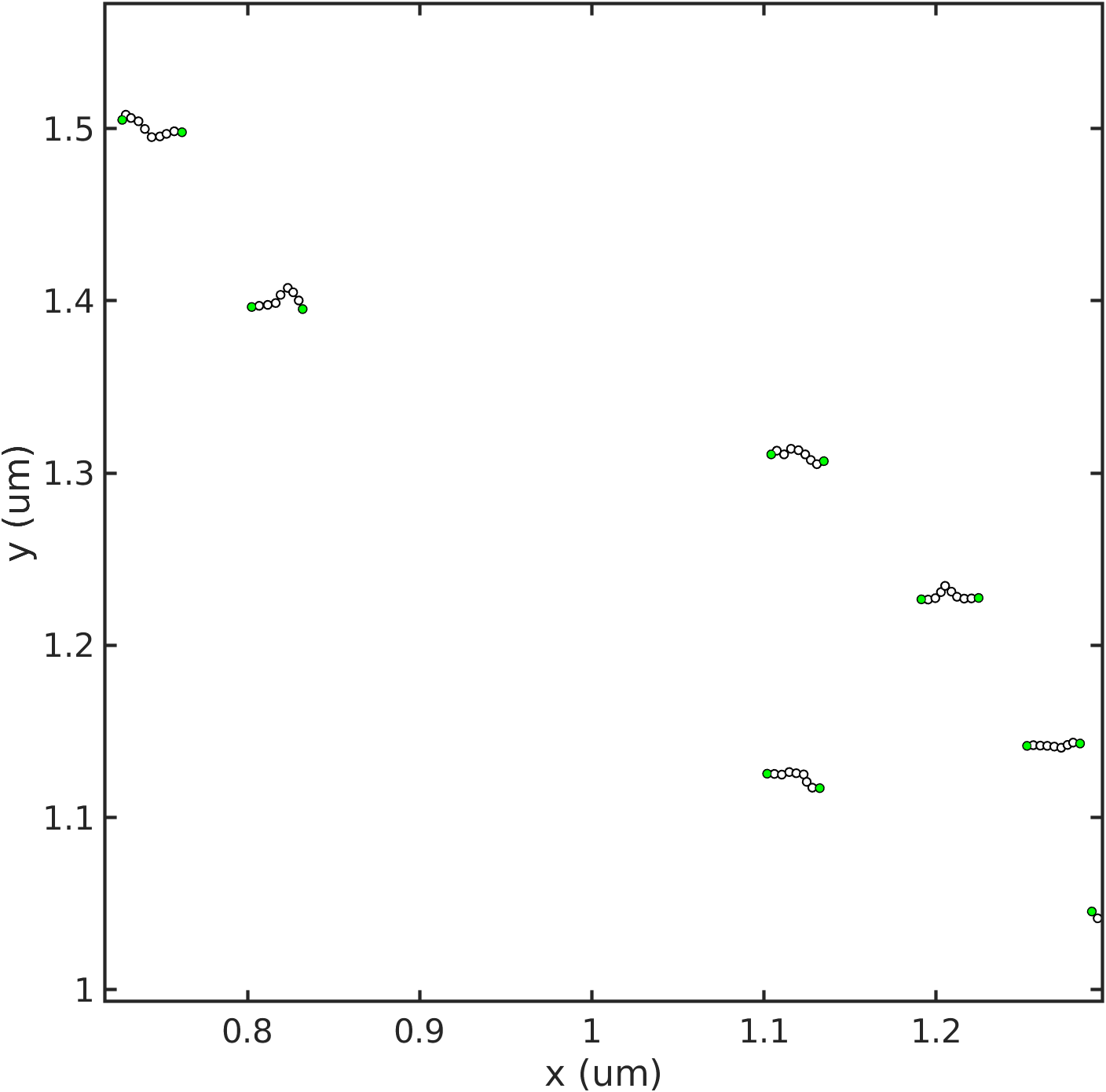}}
	\caption{The distribution diagram of the micro-crack initiation after screening} \label{fig:ms-meso-crack}
\end{figure}

\begin{figure}[h]
	\centering
	\subfigure[Iteration: 1]{\includegraphics[width=0.2\textheight]{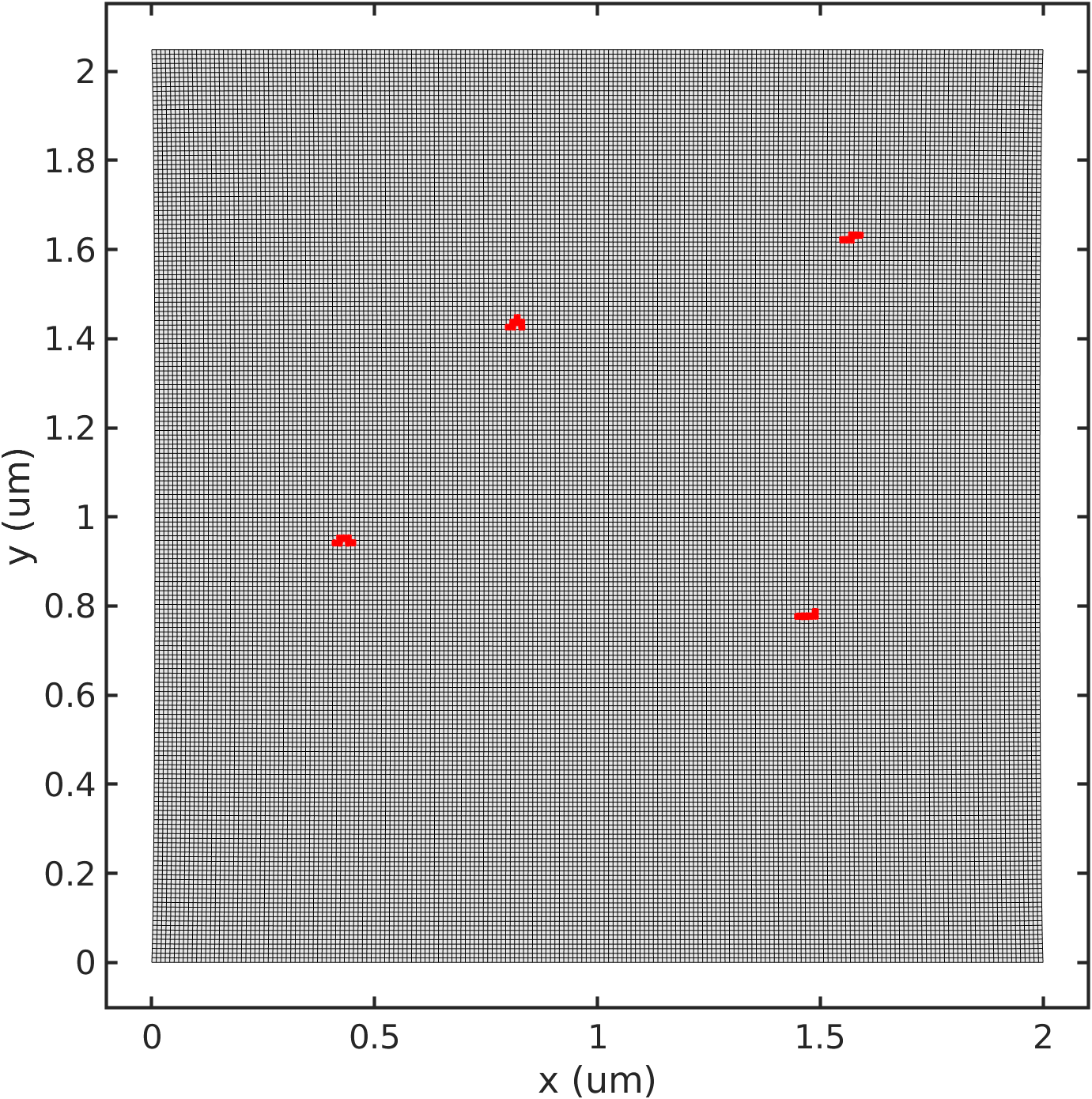}}
	\subfigure[Iteration: 5]{\includegraphics[width=0.2\textheight]{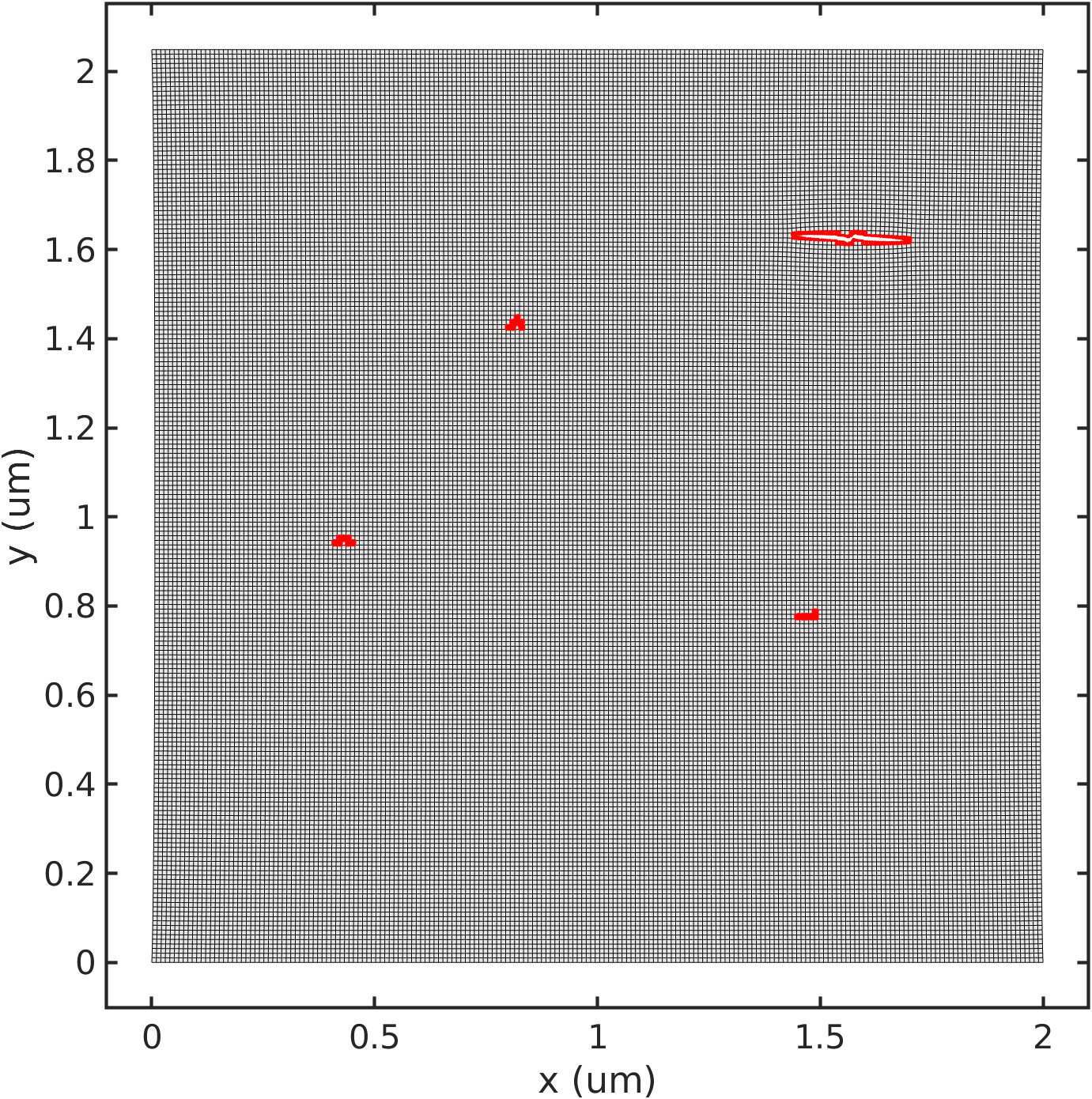}}
	\subfigure[Iteration: 15]{\includegraphics[width=0.2\textheight]{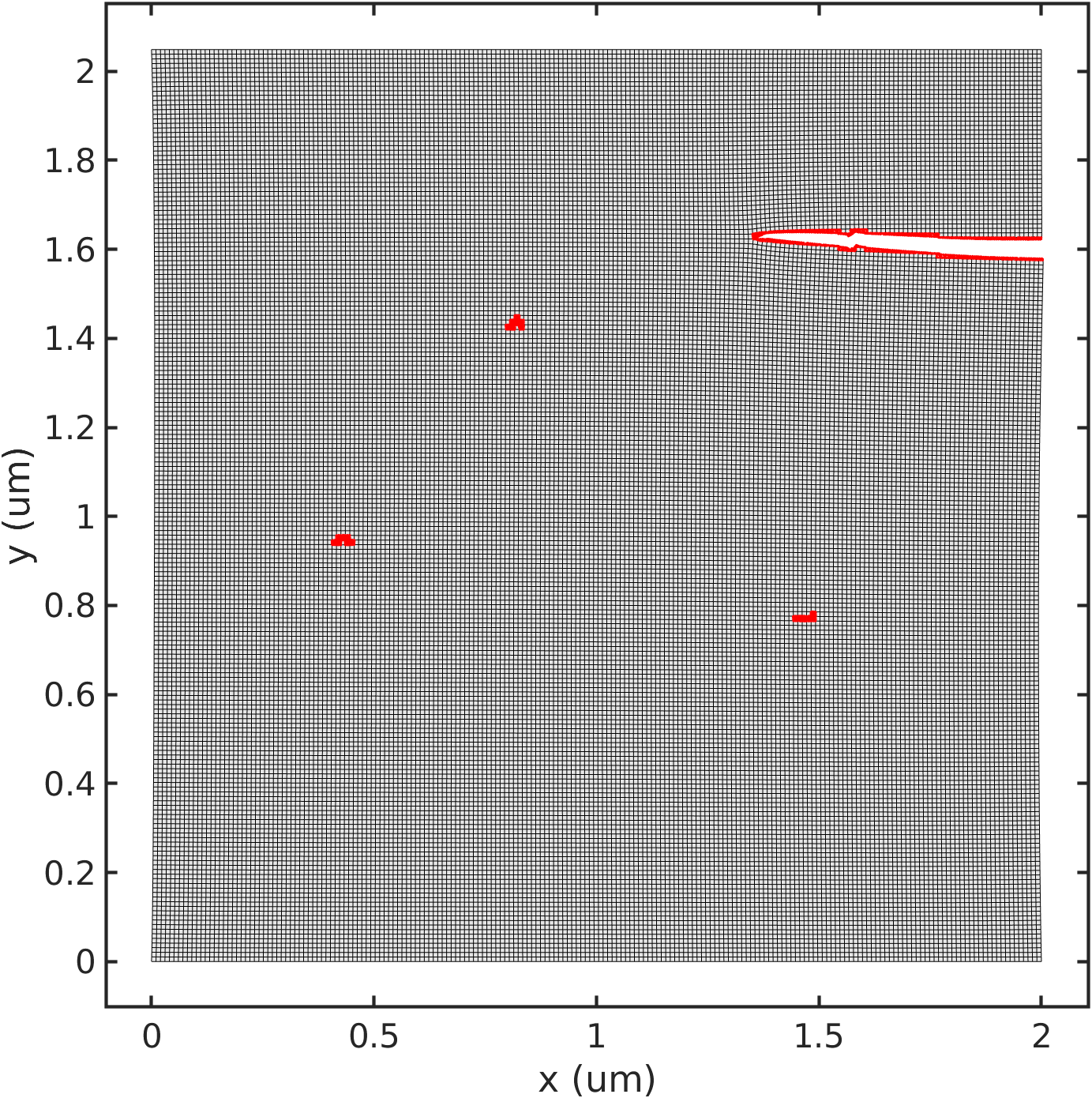}}\\
	\subfigure[Iteration: 30]{\includegraphics[width=0.2\textheight]{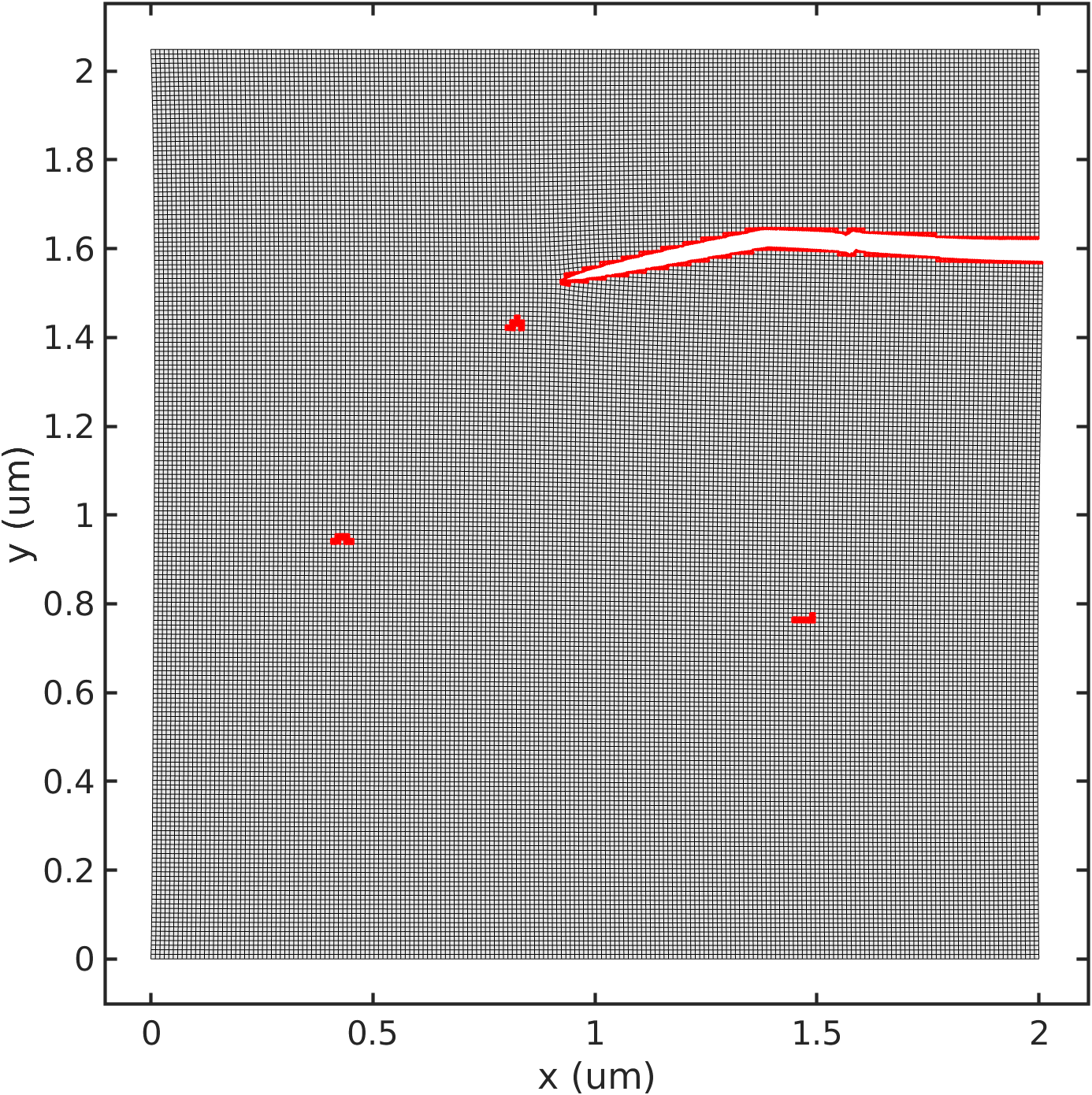}}
	\subfigure[Iteration: 45]{\includegraphics[width=0.2\textheight]{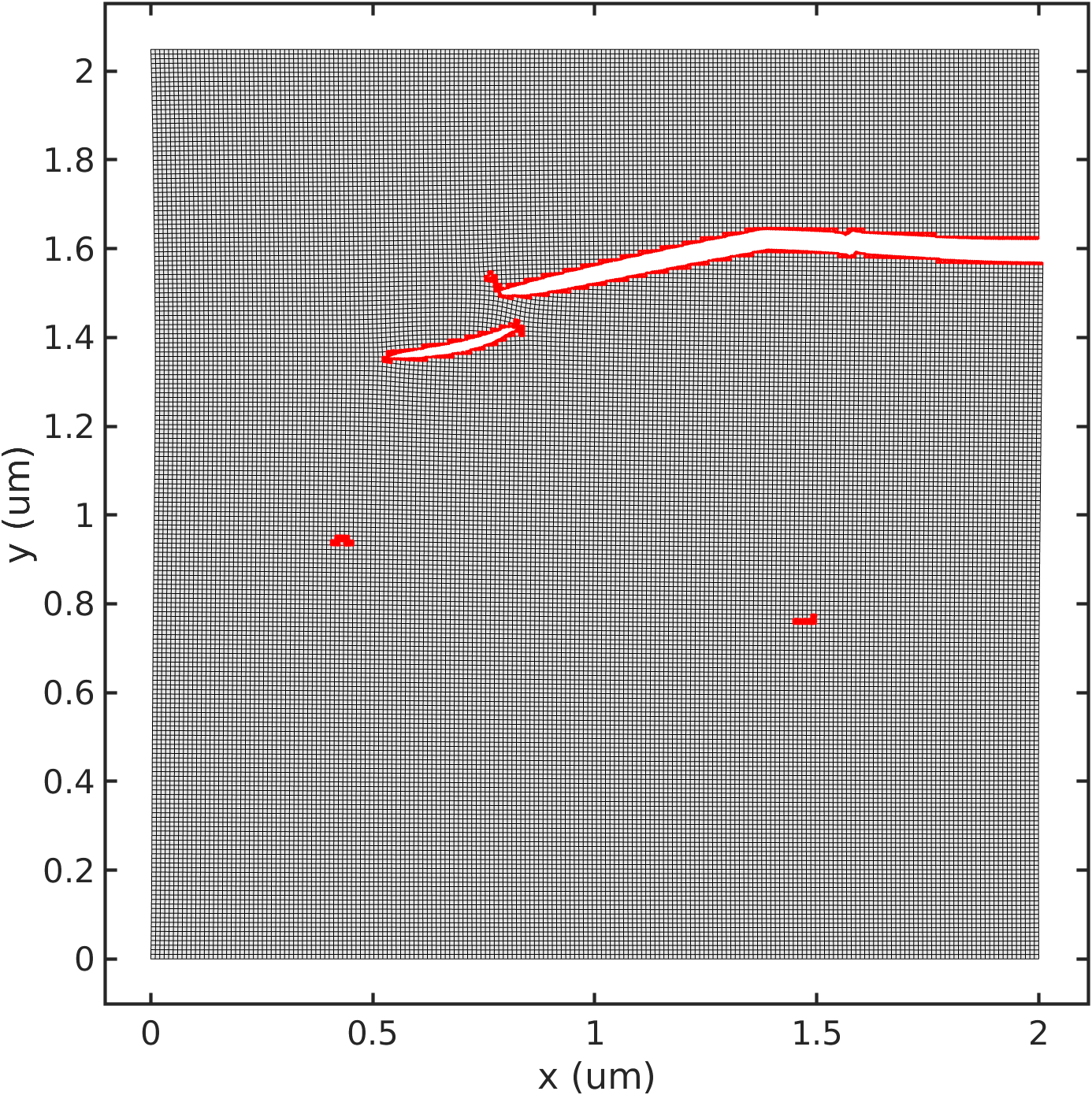}}
	\subfigure[Iteration: 65]{\includegraphics[width=0.2\textheight]{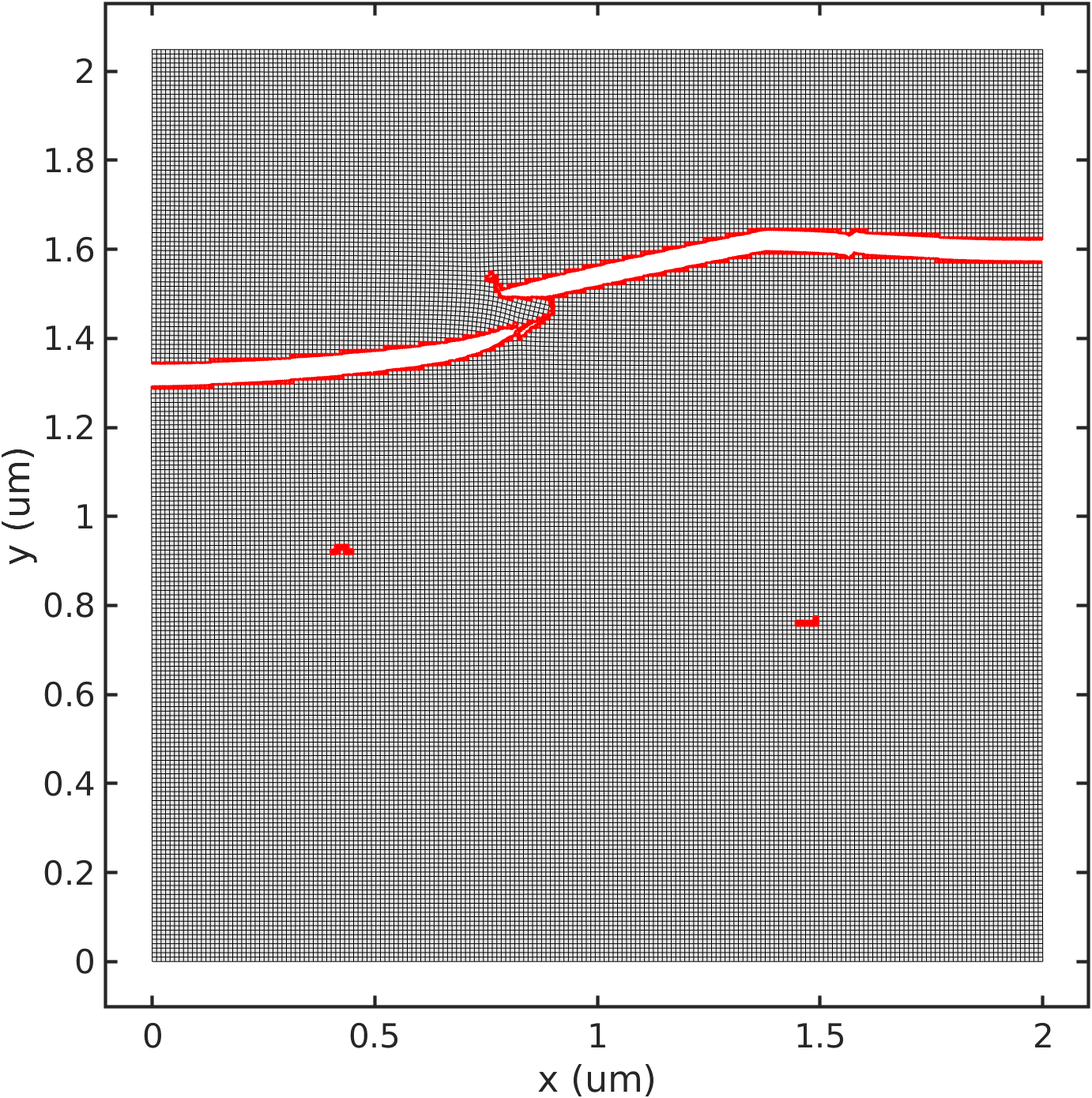}}
	\caption{The crack path of micro-crack after multi-points initiating} \label{fig:ms-micro-deform}
\end{figure}

\begin{figure}[h]
	\centering
	\subfigure[Iteration: 1]{\includegraphics[width=0.2\textheight]{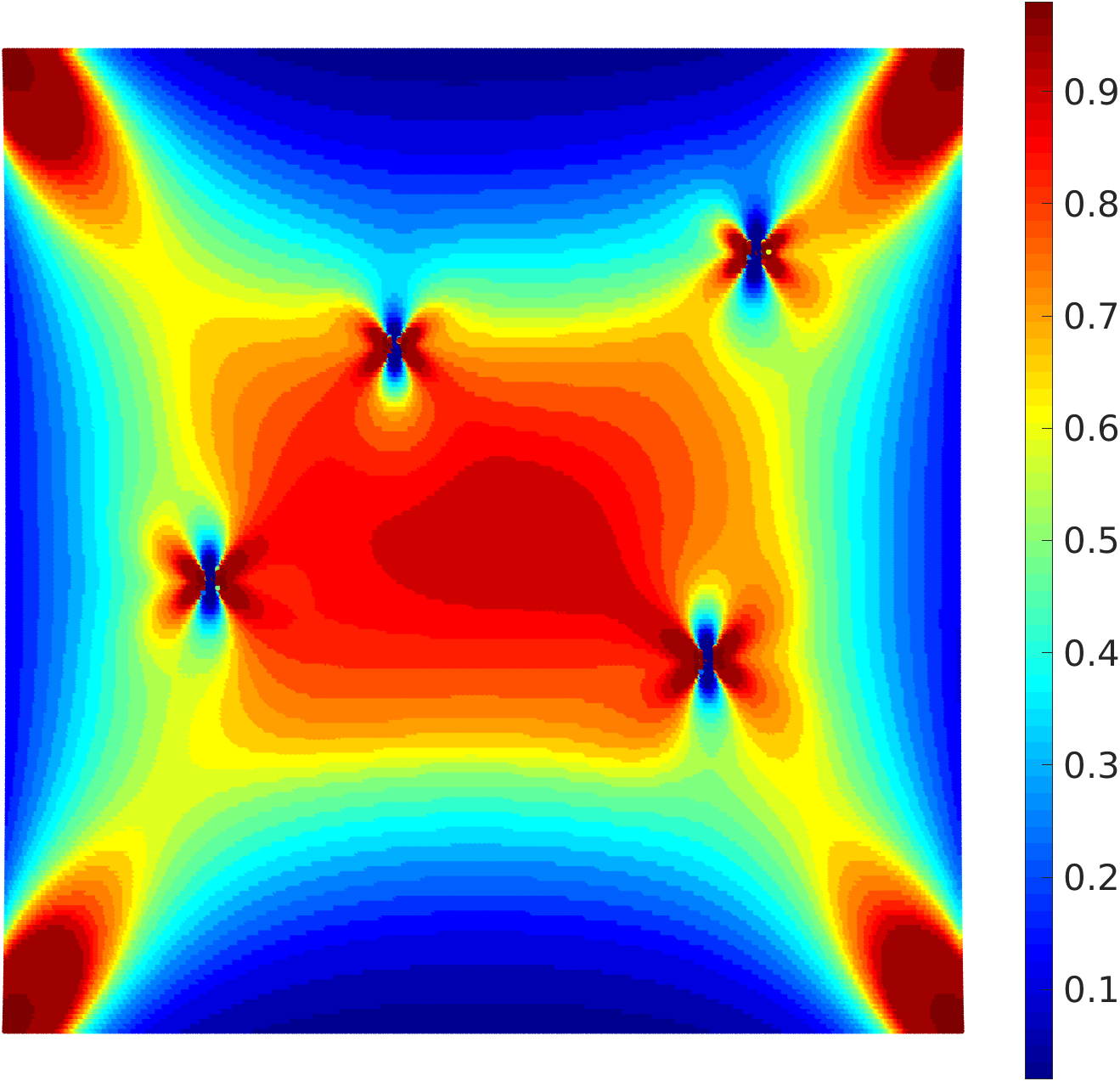}}
	\hspace{0.1in}
	\subfigure[Iteration: 30]{\includegraphics[width=0.2\textheight]{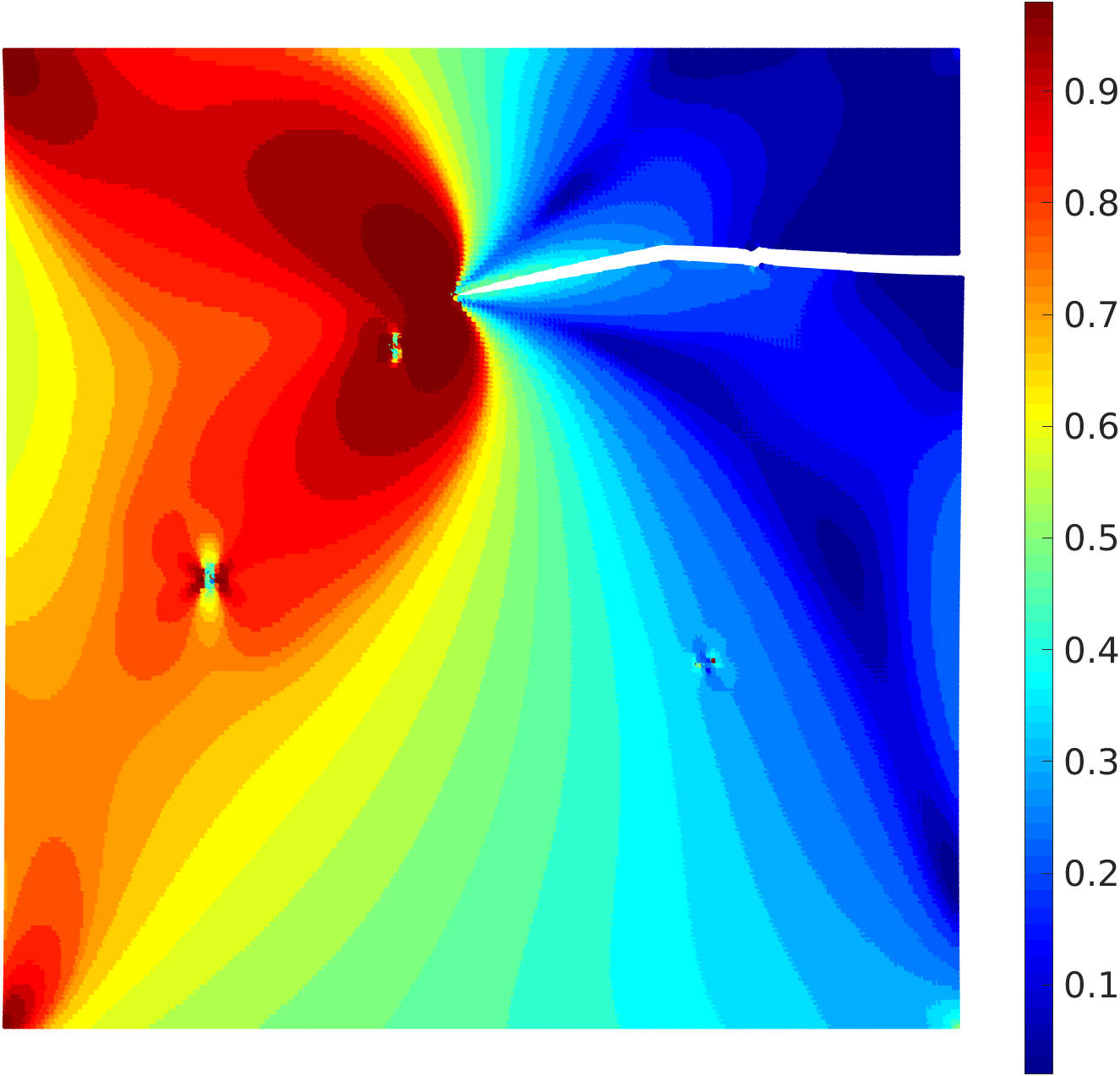}}\\
	\subfigure[Iteration: 45]{\includegraphics[width=0.2\textheight]{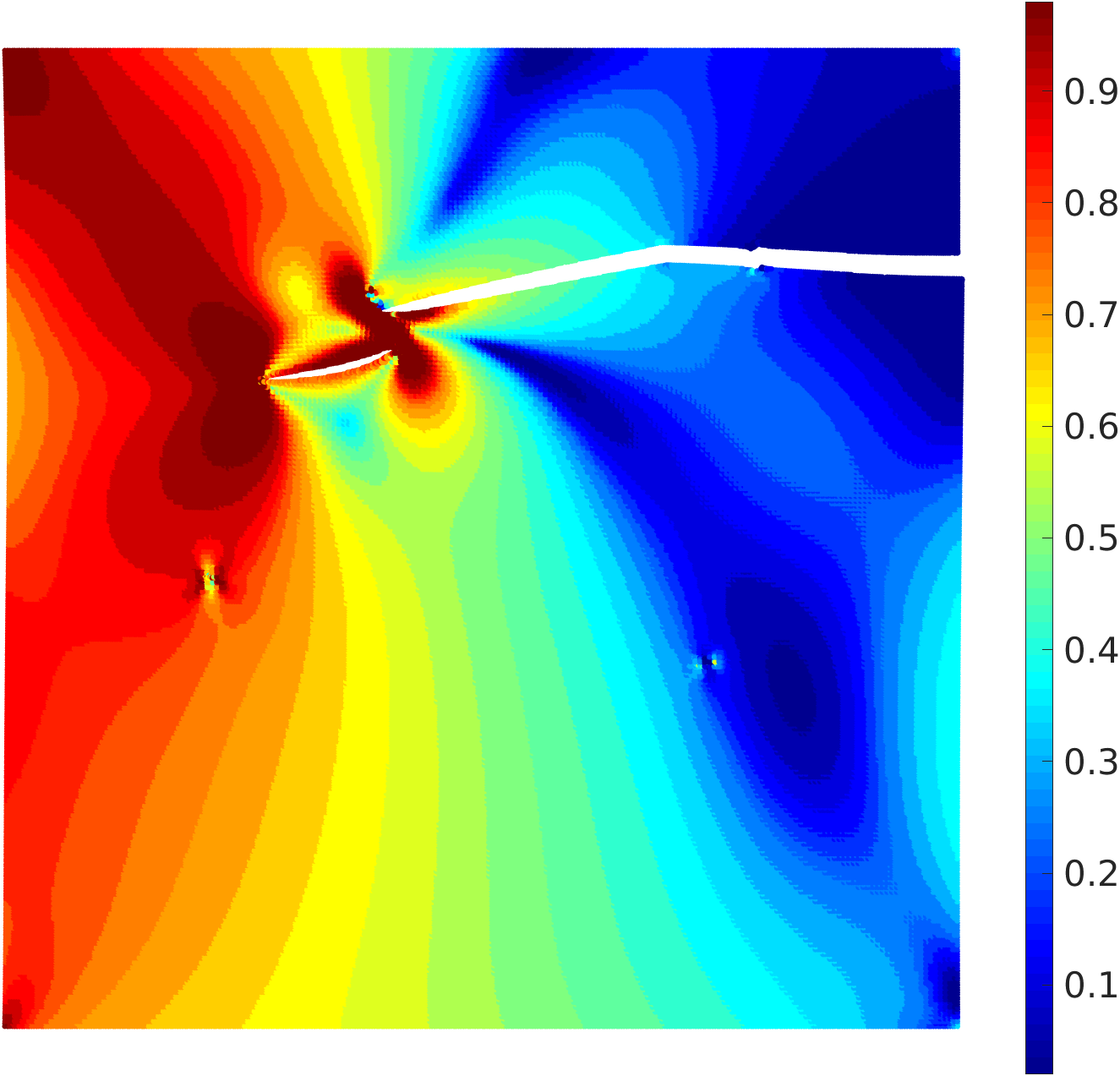}}
	\hspace{0.1in}
	\subfigure[Iteration: 65]{\includegraphics[width=0.2\textheight]{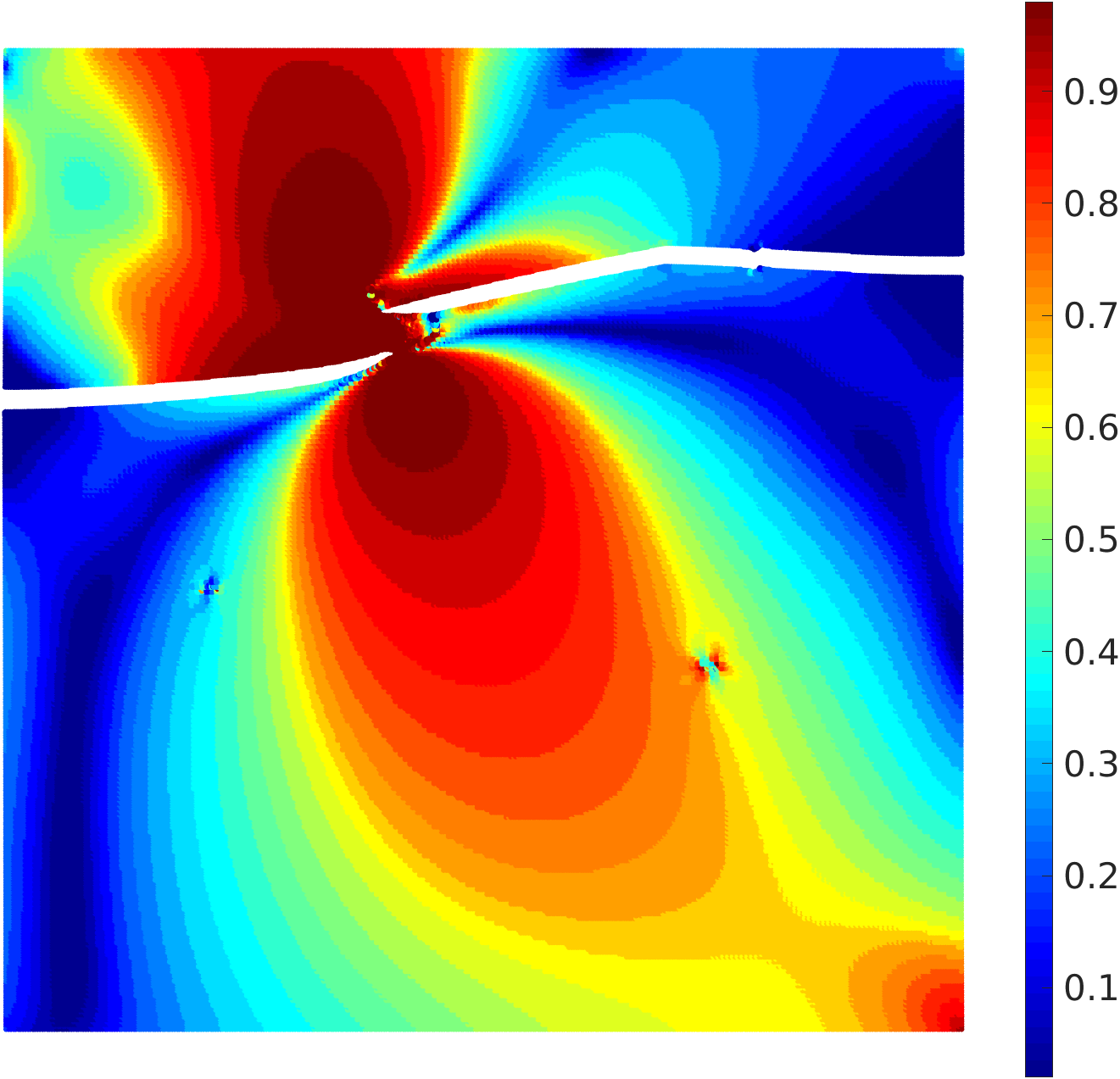}}
	\caption{The Von Mises stress contours of micro-crack after multi-point initiating} \label{fig:ms-micro-stress}
\end{figure}

At this point, the calculation of crack initiation has been completed, and the coordinate of micro-crack initiation at each sampling point has also been extracted. The next step is to analyze the propagation behavior after crack initiation. In this study, the multi-grid base XFEM is used to feedback the micro-crack initiations to the macroscopic model step by step. Firstly, we need to screen out the crack initiation from the 50 Latin hypercube sampling points  at microscopic scale. Figure \ref{fig:ms-meso-crack}  is the distribution diagram of the micro-crack initiation after screening. It can be found that only 20 samples occurred obvious crack initiation among the 50  sample points. It can be seen from the enlarged view of crack distribution, the distribution of micro-crack initiations is scattered and has different shapes. Then XFEM is used to simulate the growth process of multiple cracks at microscopic scale. Figure  \ref{fig:ms-micro-deform}  shows the crack path of micro-crack after multi-points initiating. It can be found that only  4 cracks among  20 crack initiations propagate further. At first, only one of the four main cracks grow faster. After that, another crack began to grow faster  when the main crack encountered resistance. Subsequently, those two main cracks  merged together eventually and caused a complete fracture of the crystal. Finally, the final crack is extracted as the initiation of the mesoscopic crack and transmitted to the macroscopic model. Moreover, Figure \ref{fig:ms-micro-stress}  shows the Von Mises stress contours during the crack growth process after multi-point initiating. It can be found that the global stress distribution is mainly determined by the main crack, while other micro-cracks have a certain influence on the distribution of local stress.

\begin{figure}[h]
	\centering
	\subfigure{\includegraphics[width=0.18\textheight]{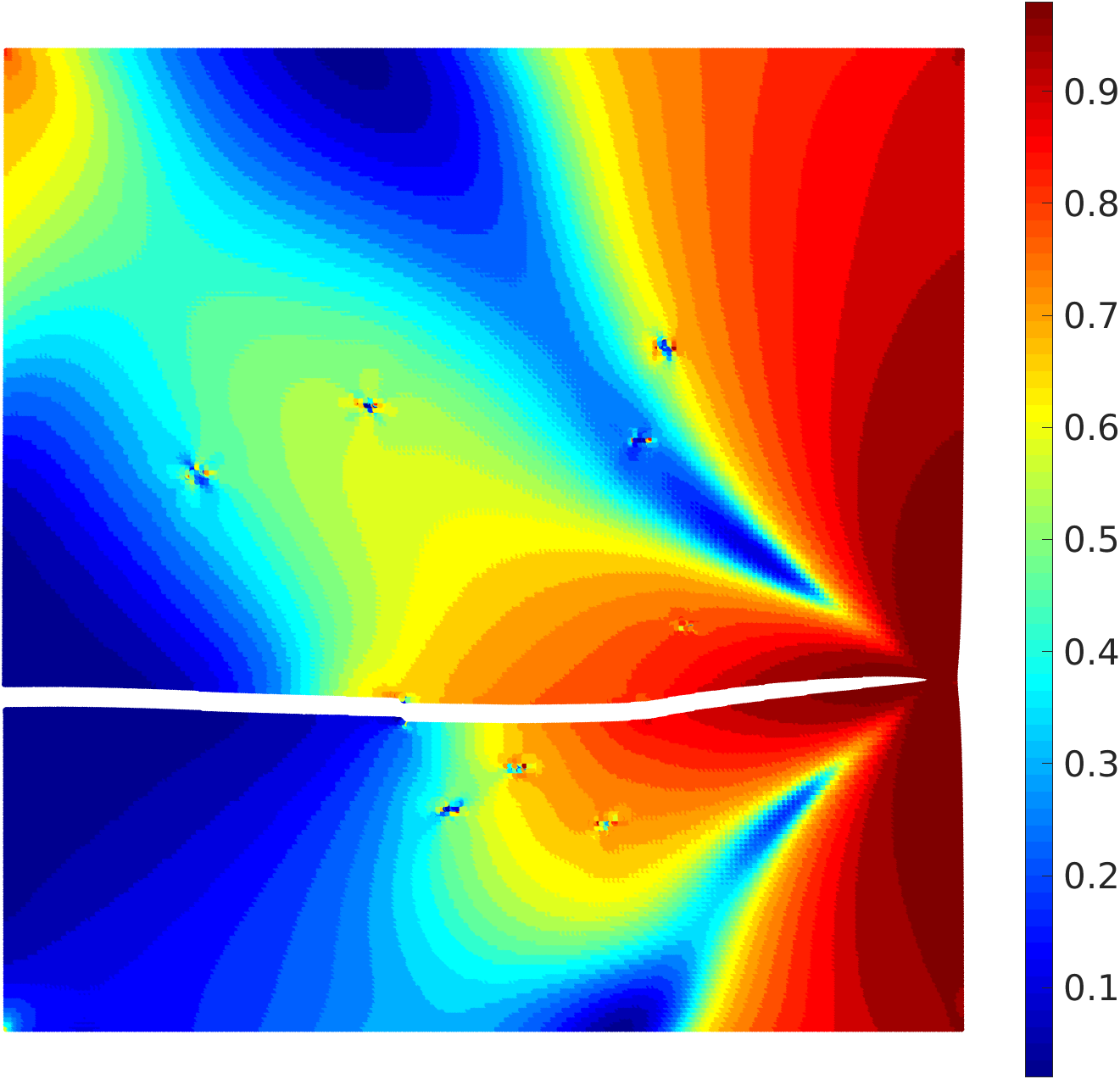}}
	\subfigure{\includegraphics[width=0.18\textheight]{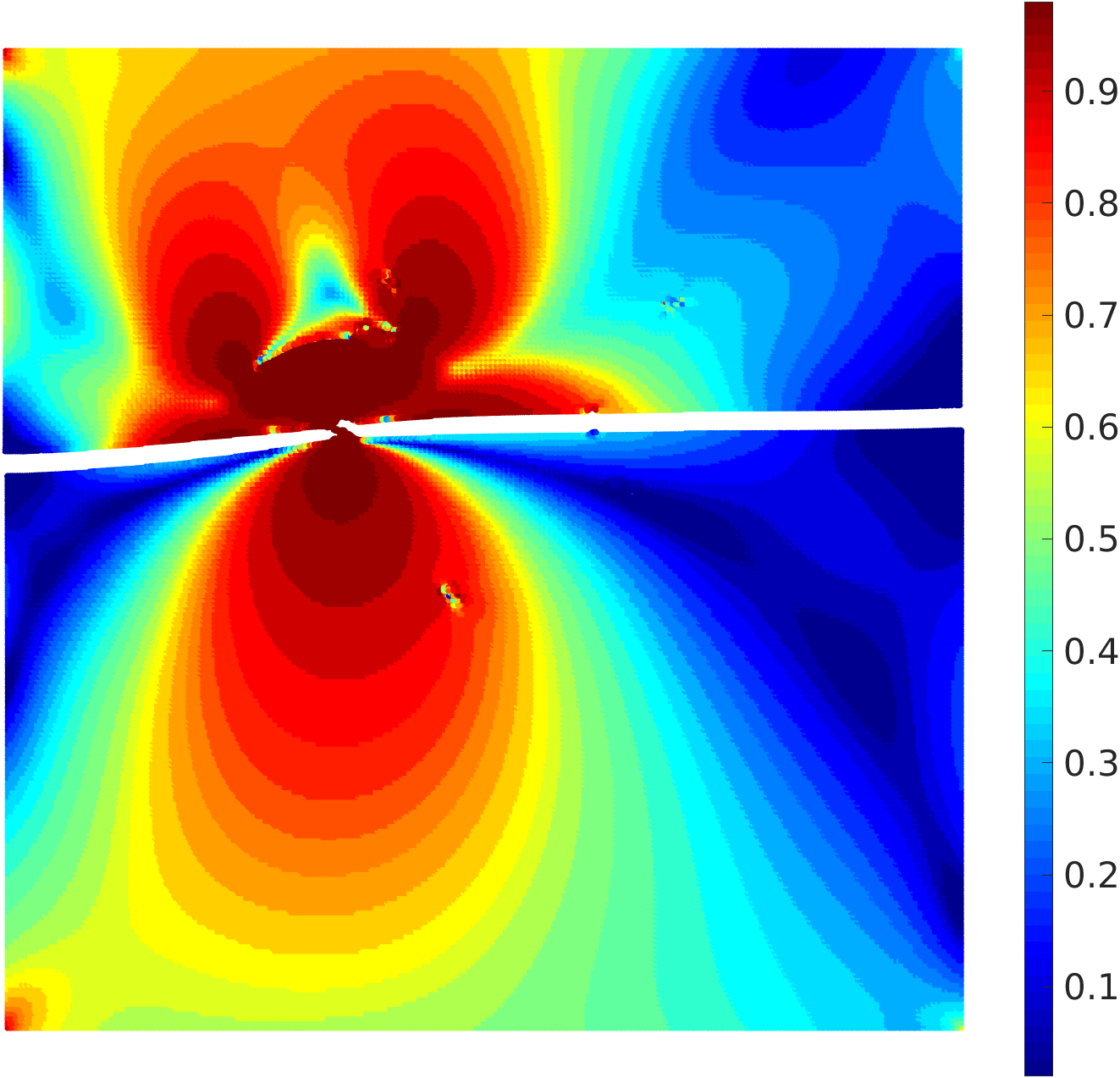}}
	\subfigure{\includegraphics[width=0.18\textheight]{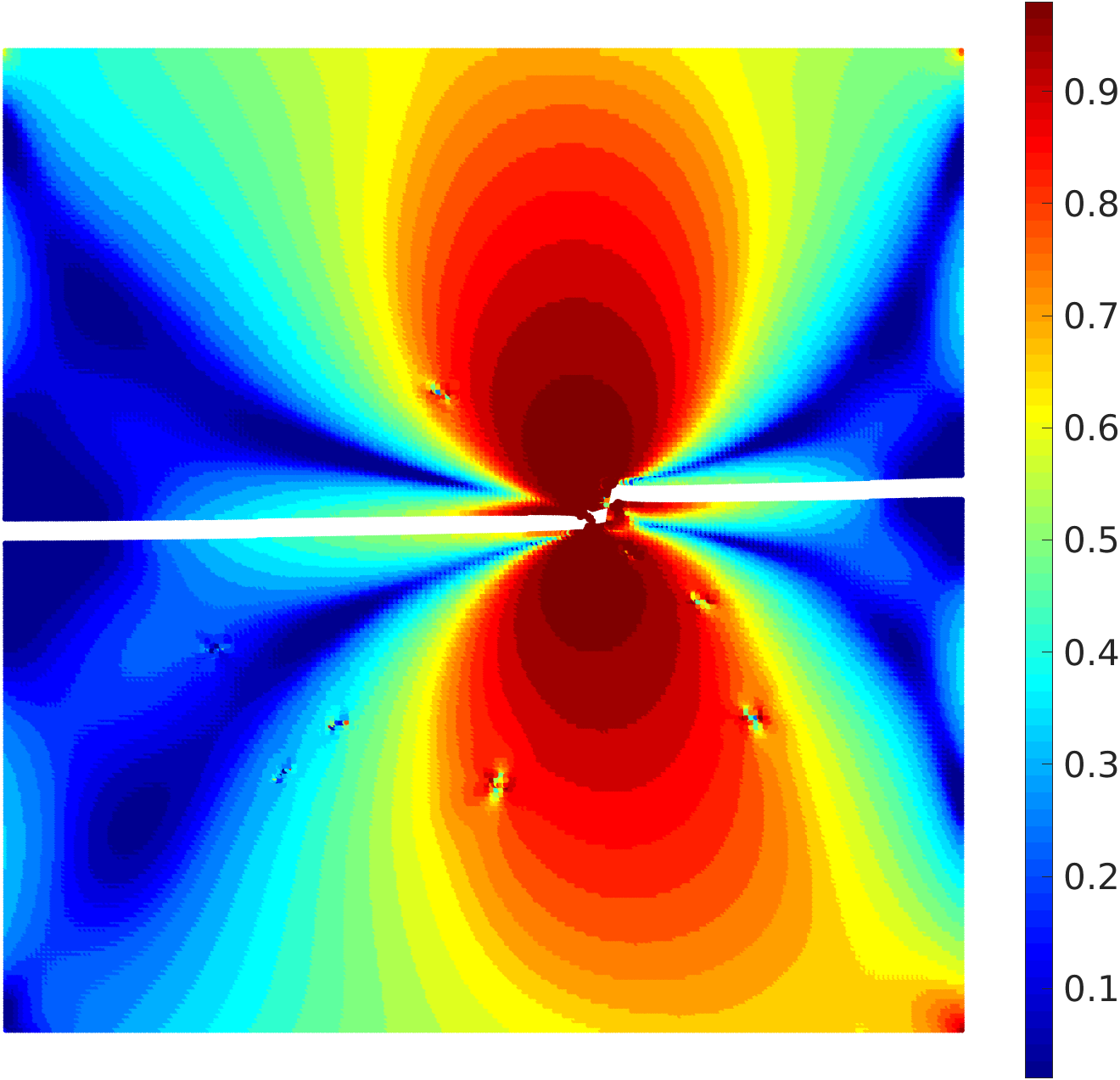}}
	\subfigure{\includegraphics[width=0.18\textheight]{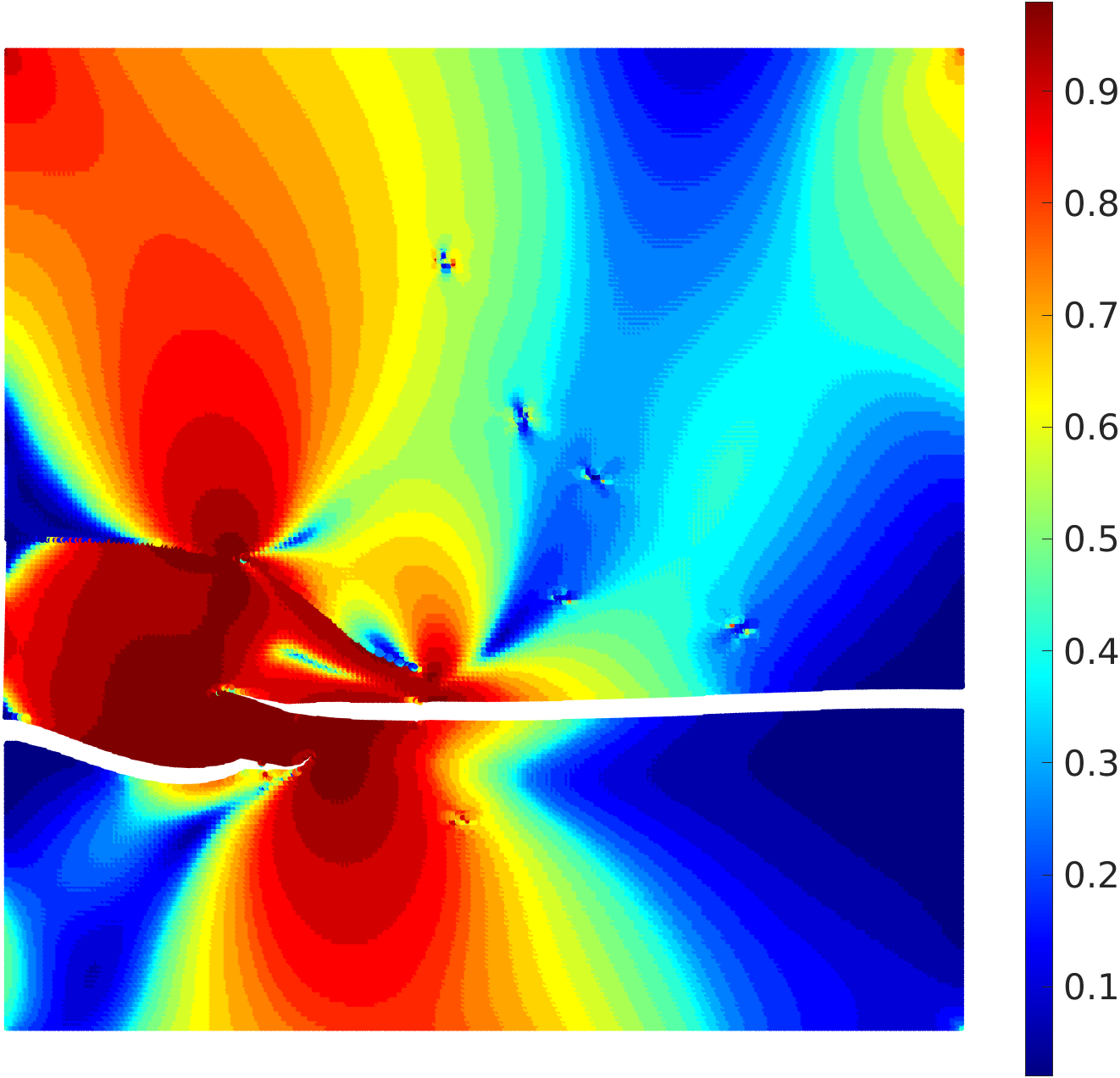}}
	\caption{The stress contours of some sample points} \label{fig:ms-meso-stress}
\end{figure}

After completing all the analysis of micro-crack initiations at 50 macro sampling points, the results of multi-crack propagation based on micro-crack initiations are obtained. The stress contours of some sample points are shown in Fig. \ref{fig:ms-meso-stress} . Then the main crack path is extracted as the mesoscopic crack initiations and transmitted to the macro model according to the final propagation path of the micro-crack under each sample point. 

\begin{figure}[h]
	\centering
	\subfigure{\includegraphics[width=0.24\textheight]{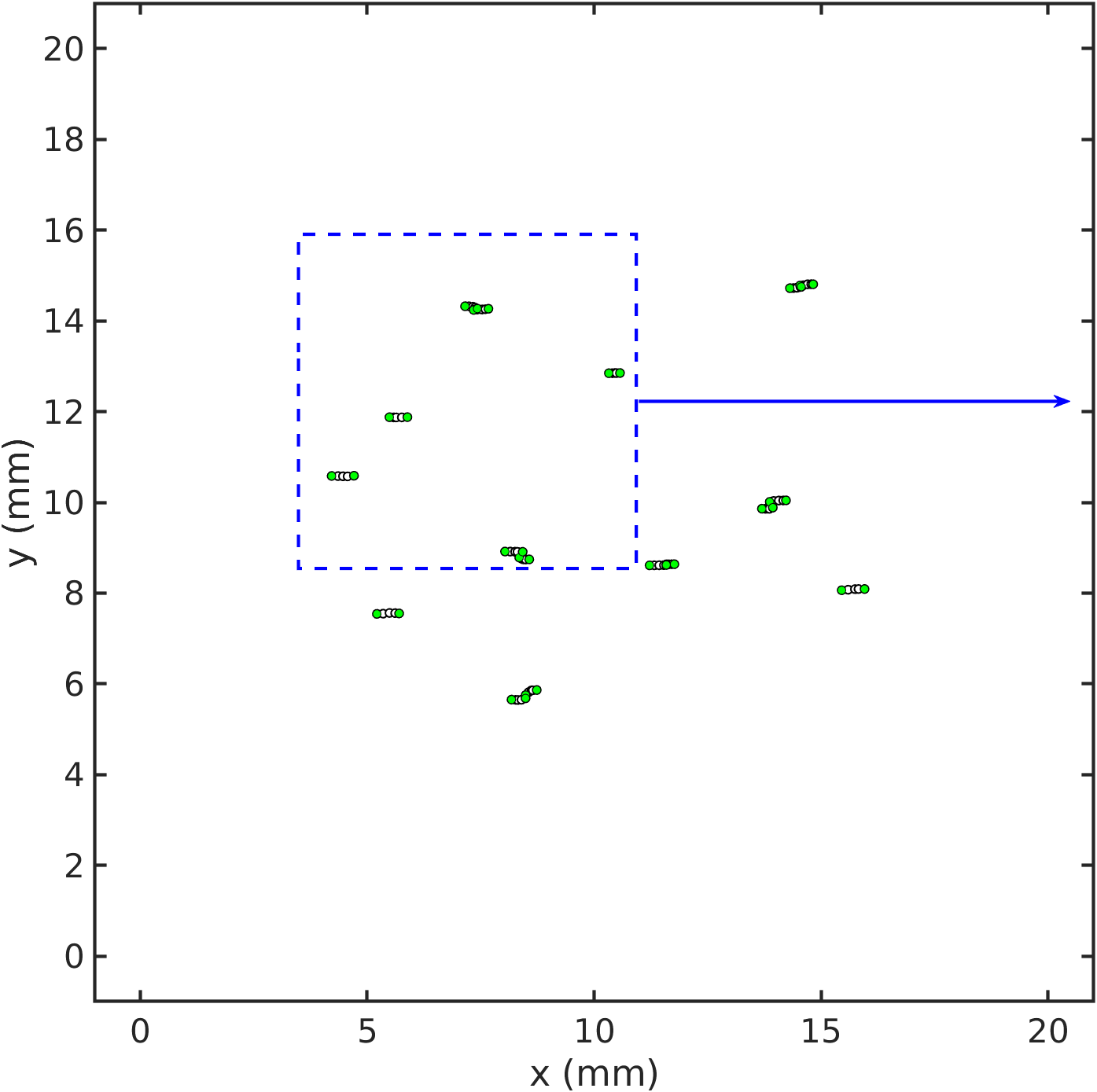}}
	\hspace{0.1in}
	\subfigure{\includegraphics[width=0.24\textheight]{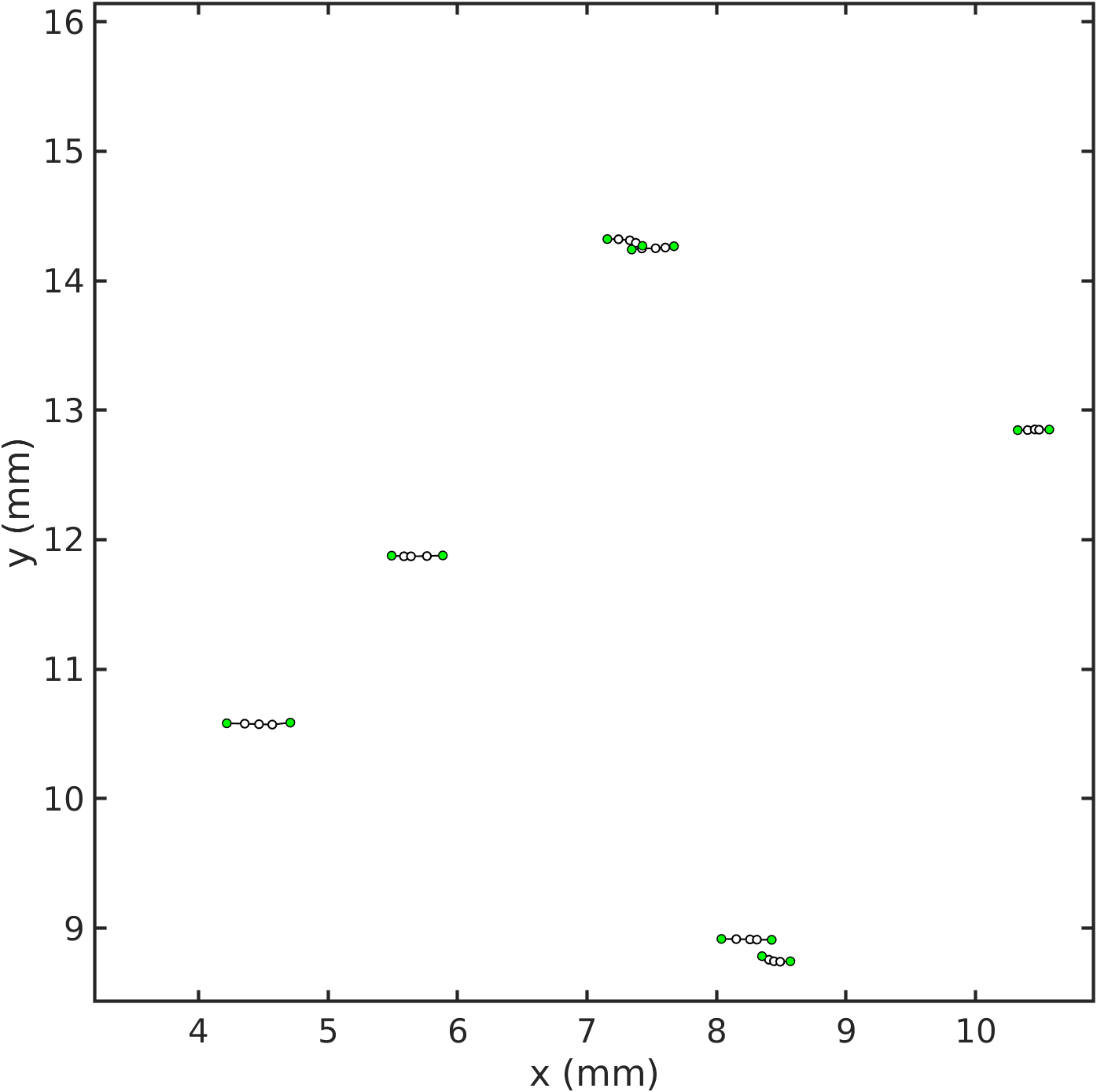}}
	\caption{The distribution of meso-crack initiations after screening} \label{fig:ms-macro-crack}
\end{figure}

\begin{figure}[h]
	\vspace{\baselineskip}
	\centering
	\subfigure[Iteration: 1]{\includegraphics[width=0.2\textheight]{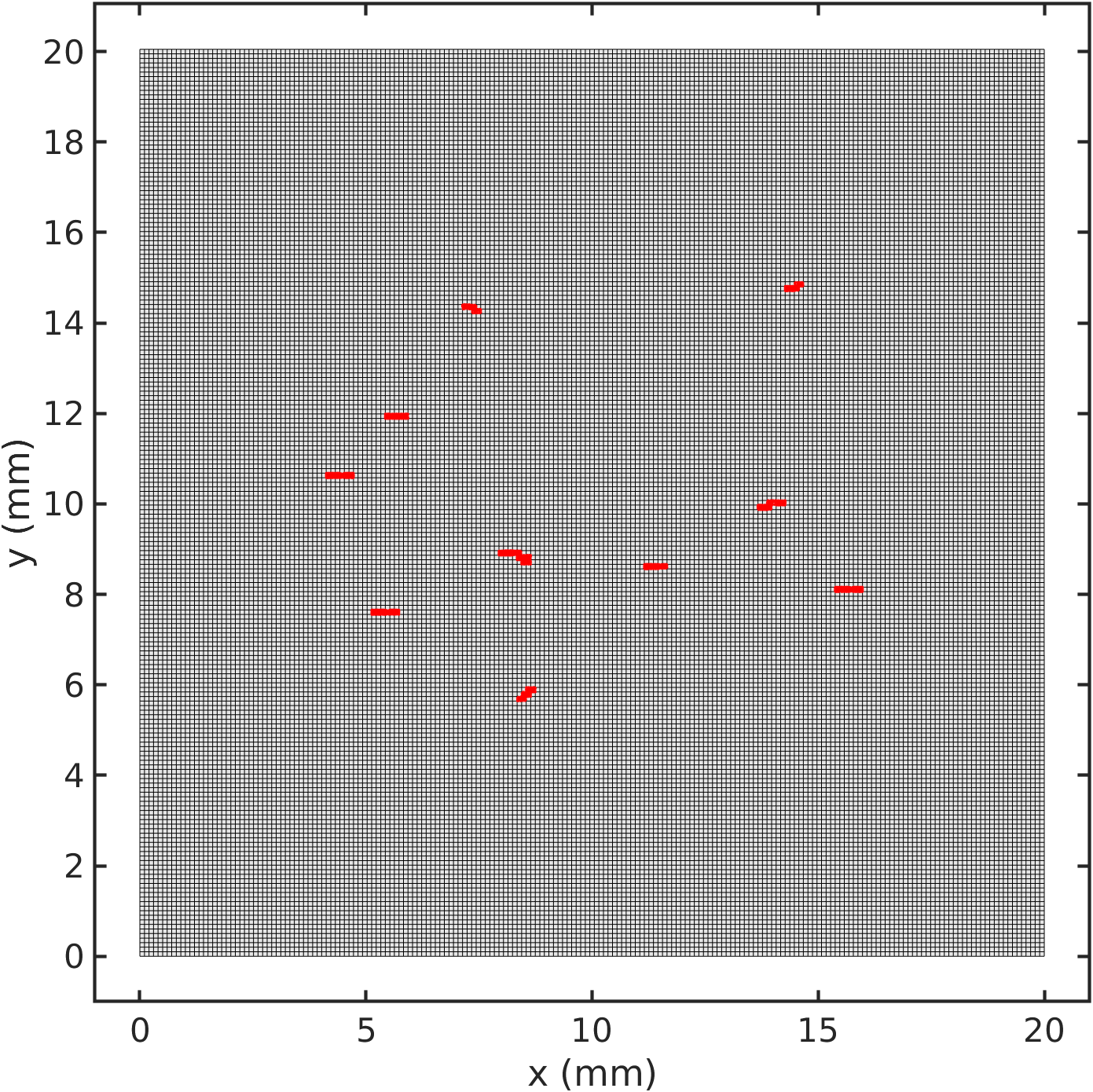}}
	\hspace{0.1in}
	\subfigure[Iteration: 10]{\includegraphics[width=0.2\textheight]{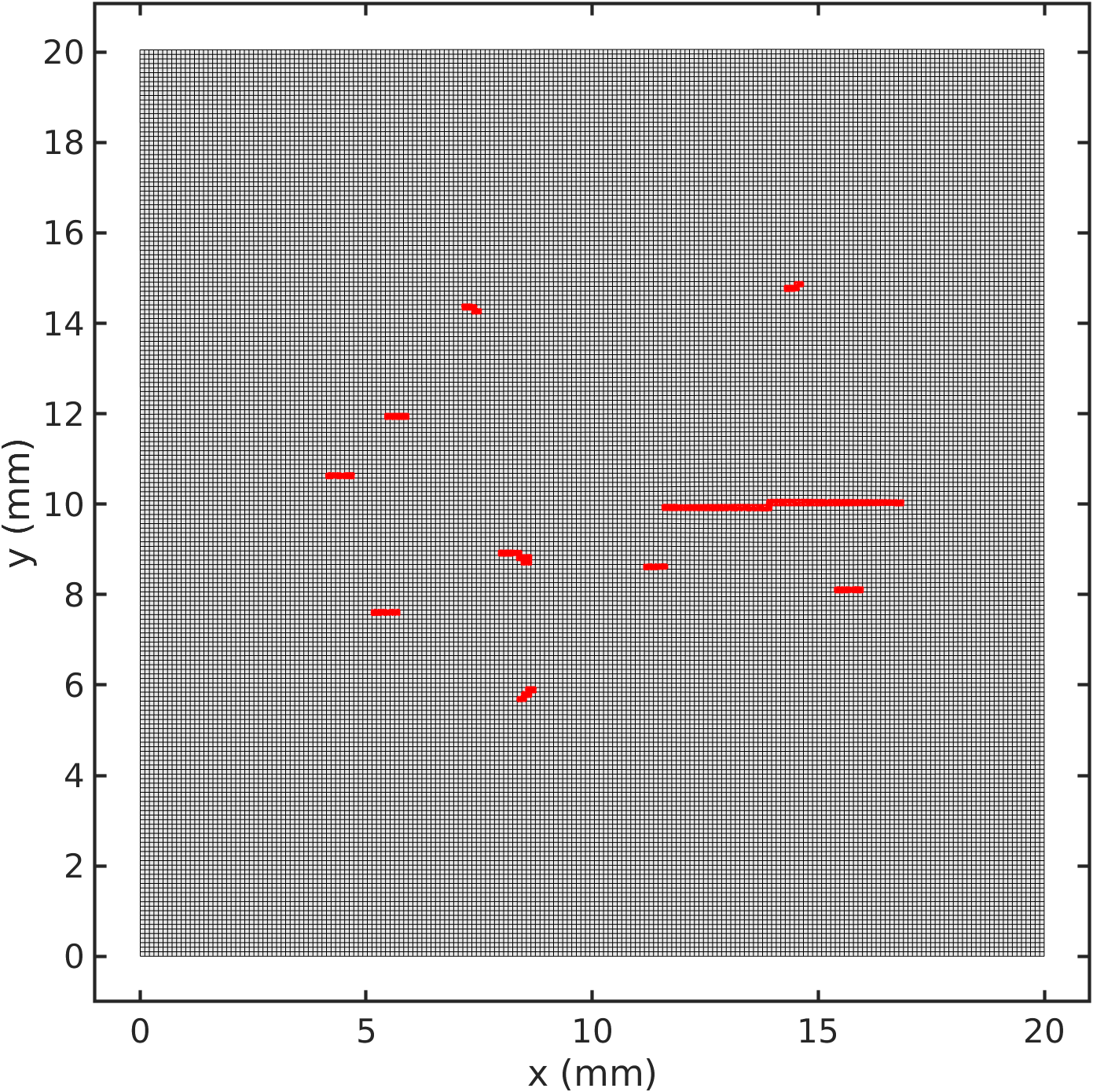}}\\
	\subfigure[Iteration: 20]{\includegraphics[width=0.2\textheight]{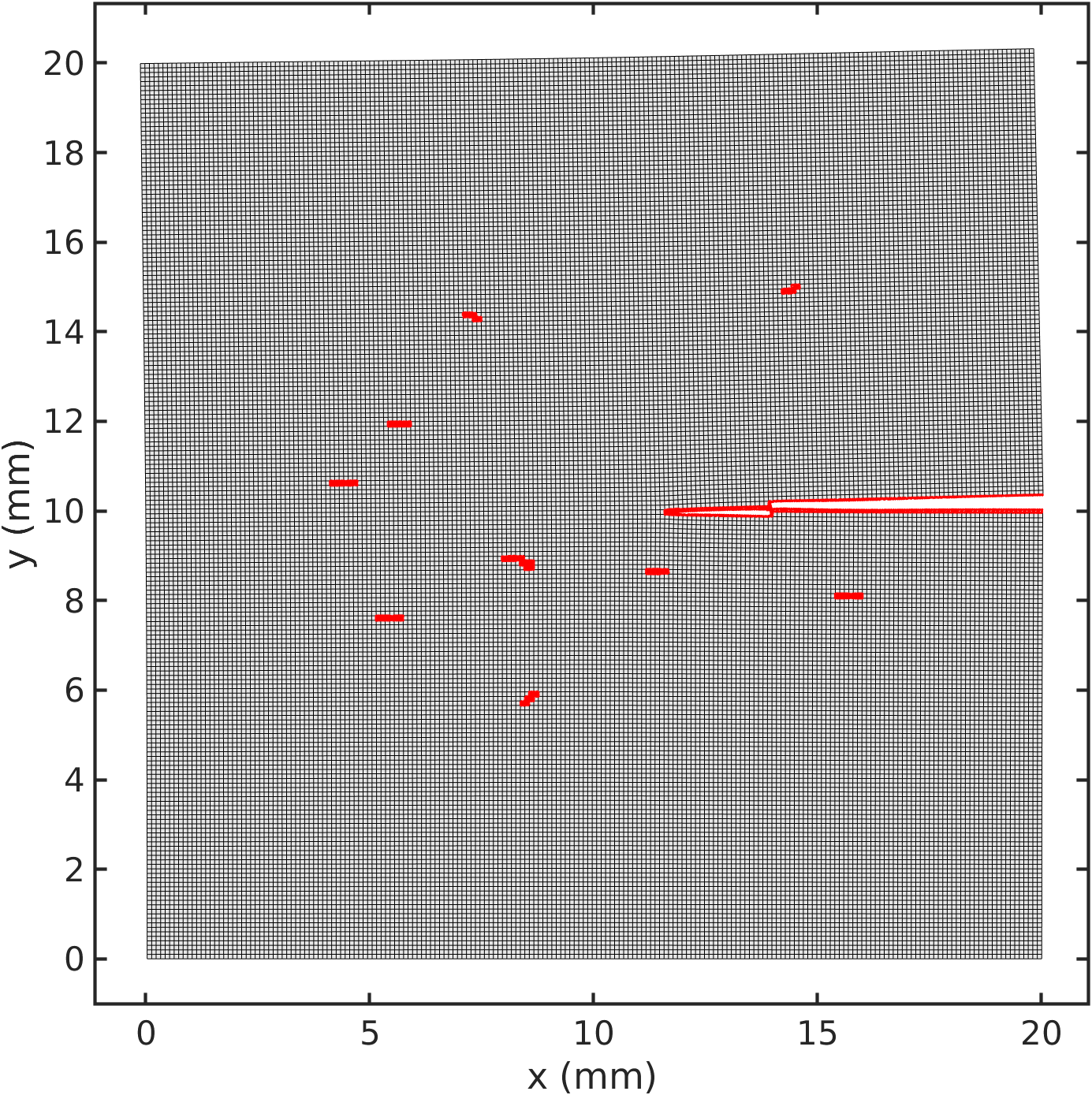}}
	\hspace{0.1in}
	\subfigure[Iteration: 30]{\includegraphics[width=0.2\textheight]{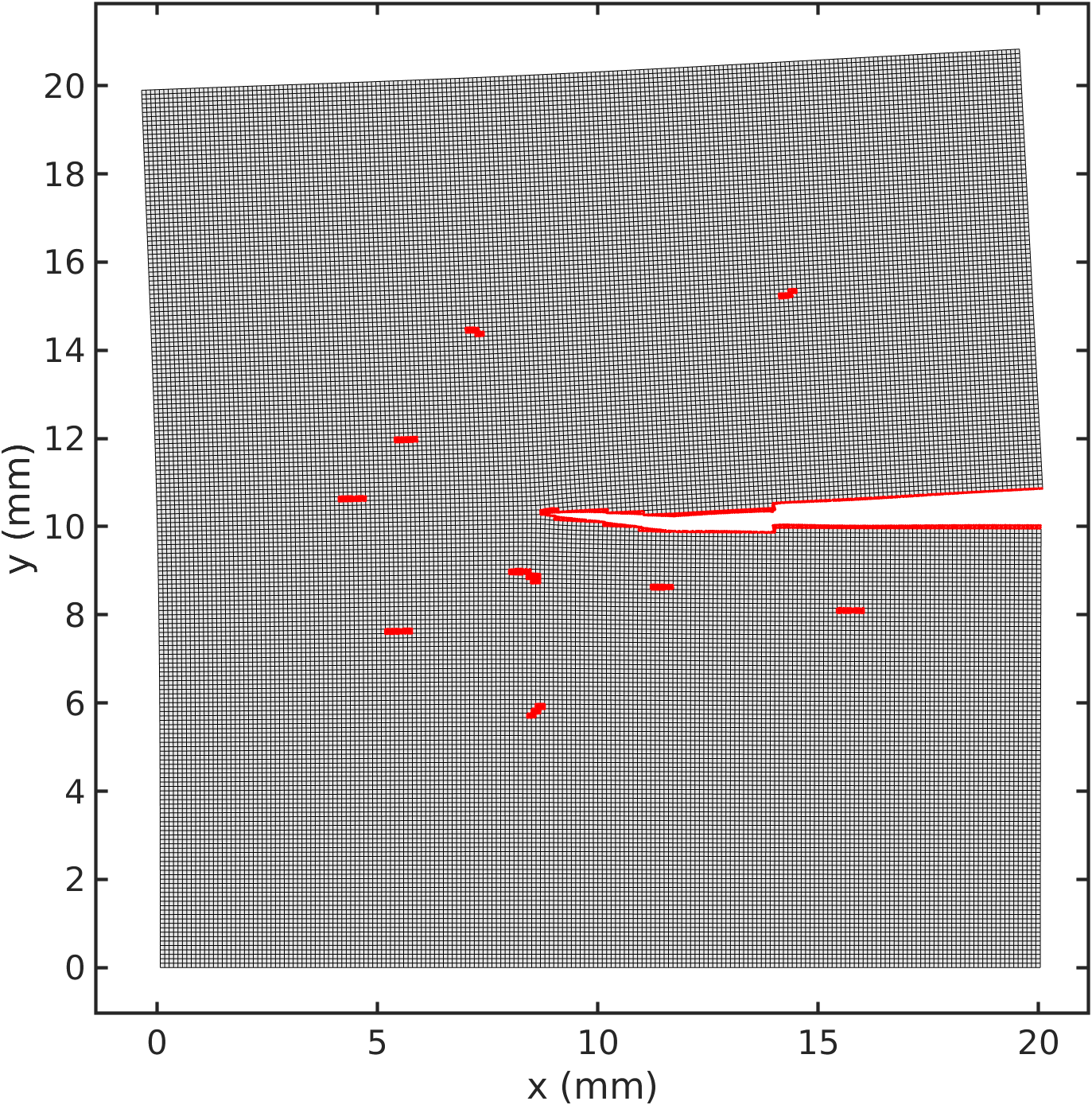}}
	\caption{The crack path of macro-crack after multi-points initiating} \label{fig:ms-macro-deform}
	\vspace{\baselineskip}
\end{figure}

Subsequently, the macro extended finite element model was constructed according to the coordinates of mesoscopic crack initiations of  50 sampling points. Then, the final macro crack propagation path was calculated by the XFEM. Figure \ref{fig:ms-macro-crack}  is the distribution of meso-crack initiations after screening. Obviously, there are only 11 sampling points with obvious crack initiations among the 50  sampling points .  Figure \ref{fig:ms-macro-deform}  shows the process of multi-crack propagation at macro scale. Figure \ref{fig:ms-macro-deform}(a), (b), (c)  and  (d)  are the crack propagation paths at the 1st, 10th, 20th and 30th iterations, respectively.

\begin{figure}[h]
	\vspace{\baselineskip}
	\centering
	\subfigure[Iteration: 1]{\includegraphics[width=0.2\textheight]{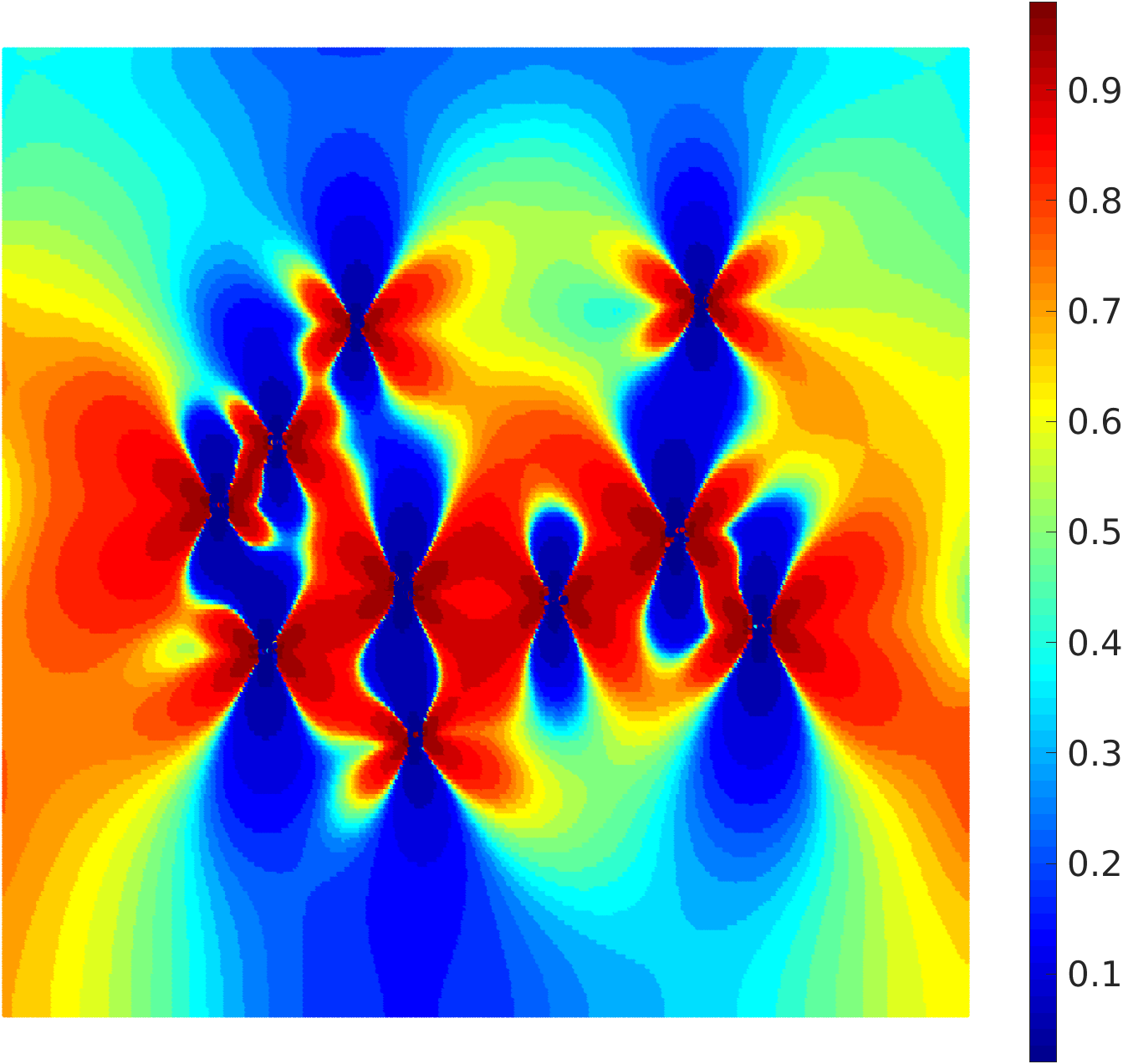}}
	\hspace{0.1in}
	\subfigure[Iteration: 10]{\includegraphics[width=0.2\textheight]{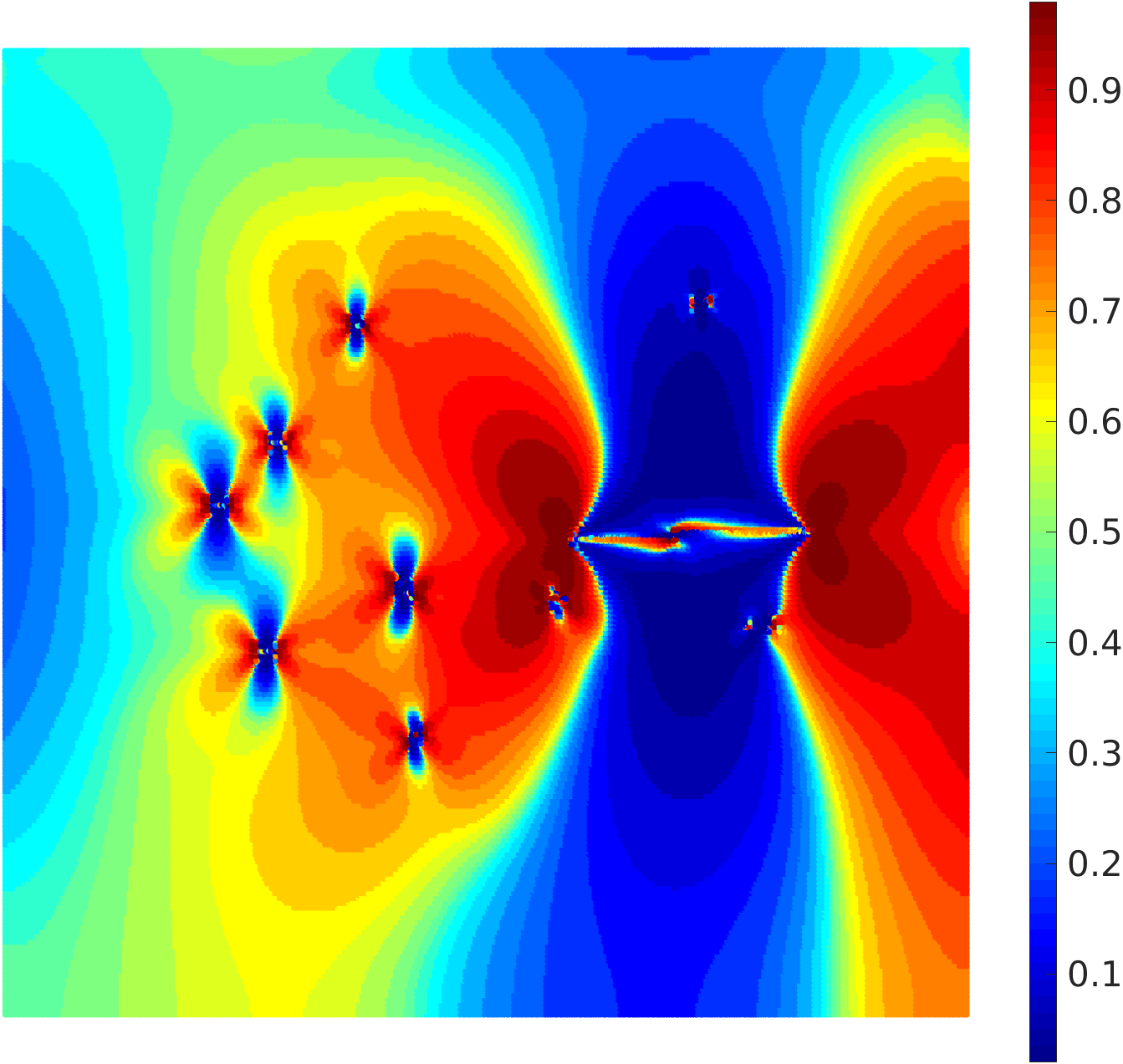}}\\
	\subfigure[Iteration: 20]{\includegraphics[width=0.2\textheight]{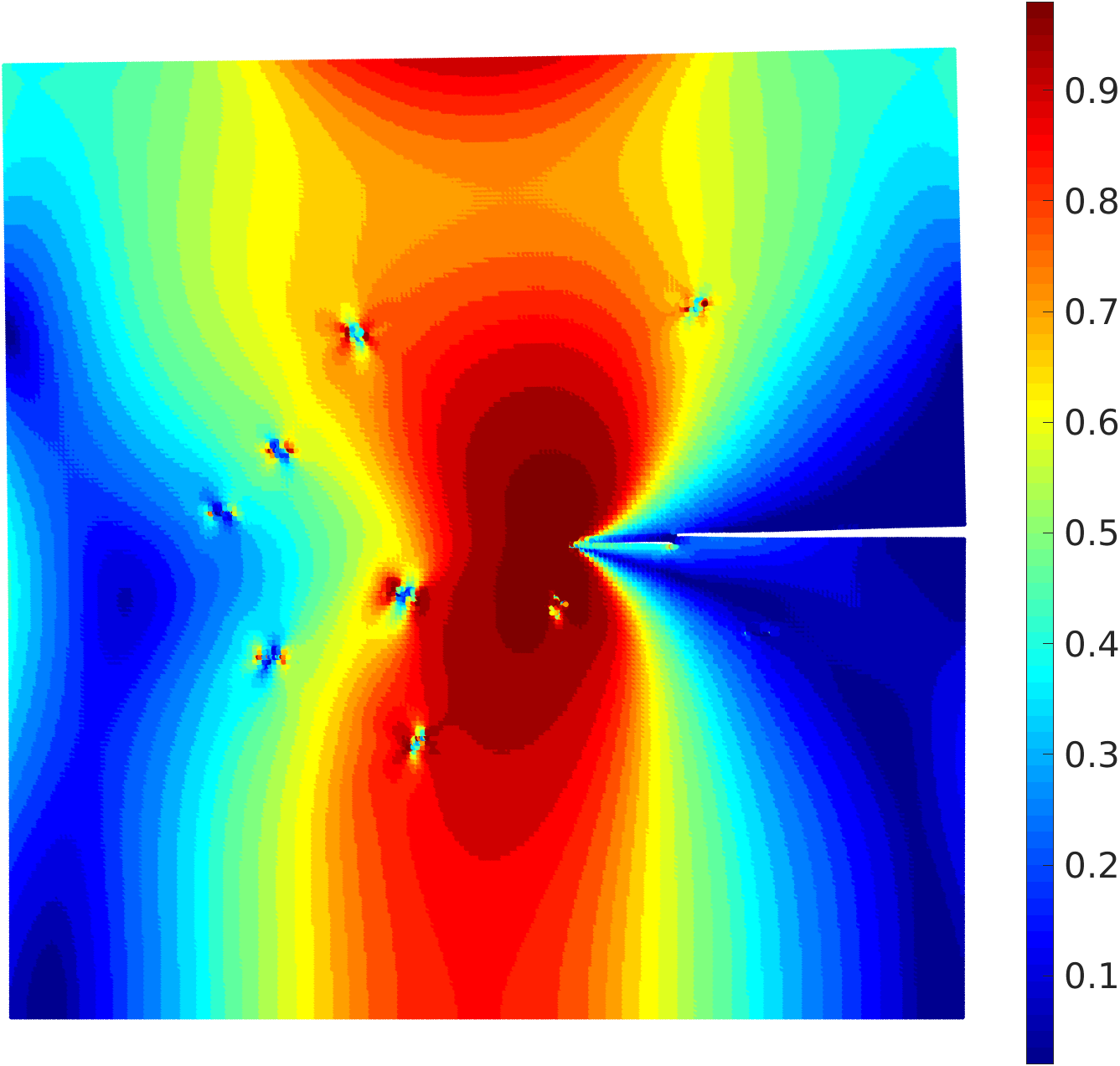}}
	\hspace{0.1in}
	\subfigure[Iteration: 1]{\includegraphics[width=0.2\textheight]{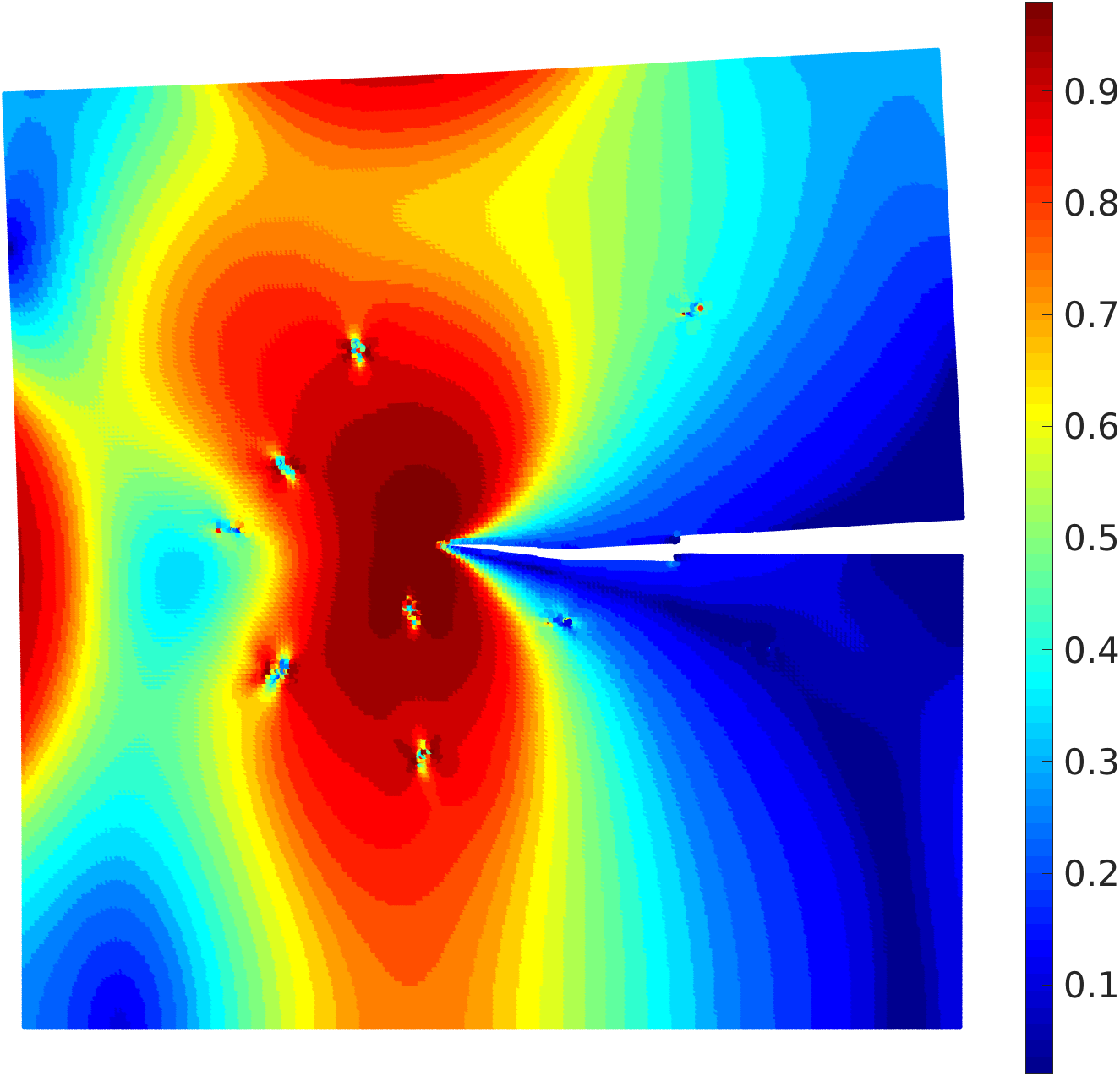}}
	\caption{The Von Mises stress contours of macro-crack after multi-point initiating} \label{fig:ms-macro-stress}
	\vspace{\baselineskip}
\end{figure}

It can be found that 10 of 11crack initiations propagate further. Moreover, only one of the cracks grows faster. When the main crack arrives at the edge of the plate, the plate appears obvious deformation.The Von Mises stress contours during the process of macro crack propagation after multi-point initiating are shown in Fig. \ref{fig:ms-macro-stress}.  It can be found that,  the stress concentration appeared at all the crack initiation locations at the initial stage. Then the stress concentration of the main crack increased significantly and covered the stress distribution at other crack initiations with the loading increasing, but the stress concentration phenomenon is still existed at each crack initiations in the local region.

\section{conclusions}
In this study, a multi-grid sampling multi-scale method is proposed to analyze the crack propagation. Cracks are analyzed comprehensively on both macroscopic and microscopic scales. The influence of micro-structure is fully considered while ensuring the size of macroscopic model. The work is summarized as follows:
	\begin{enumerate}[(1)]
		\item A multi-grid sampling multi-scale method is proposed by combining with finite element, extended finite element and molecular dynamics methods.The process of crack propagation is analyzed comprehensively from the microscopic initiation to the macroscopic propagation.		
		\item Multi-grid finite element method is used to transmit the macroscopic displacement boundary conditions to the atomic model and the multi-grid extended finite element method is used to feedback the microscopic crack initiations to the macroscopic model.		
		\item An image recognition based crack extracting method is proposed to extract the crack coordinate from the LAMMPS result files efficiently.It can not only save a lot of storage space, but also largely reduces the computational cost of crack extracting process.
		\end{enumerate}
	
	\section*{acknowledgements}
	This work was partially supported by Project of the Key Program of National Natural Science Foundation of China (Grant Number 11972155).
	
	\bibliographystyle{unsrt}
	\bibliography{bibfile}

\end{document}